\RequirePackage{fix-cm}
\documentclass{svjour3}   
\smartqed  

\newtheorem{algorithm}{Algorithm}
\def\Frac#1#2{\frac{\displaystyle{#1}}{\displaystyle{#2}}}
\usepackage{graphicx}

\begin{document}

\title{Fast, reliable and unrestricted iterative computation of Gauss--Hermite and Gauss--Laguerre quadratures\thanks{This work was supported by  {\emph{Ministerio de Ciencia, Innovaci\'on y Universidades}}, projects MTM2015-67142-P (MINECO/FEDER, UE) and PGC2018-098279-B-I00 (MCIU/AEI/FEDER,UE)
}}

\titlerunning{Iterative computation of Gauss--Hermite and Gauss--Laguerre quadratures} 

\author{Amparo Gil \and Javier Segura \and Nico~M.~Temme}
\authorrunning{A. Gil, J. Segura, N. M.~Temme} 

\institute{Amparo Gil \at
              Departamento de Matem\'atica Aplicada y CC. de la Computaci\'on.
ETSI Caminos. Universidad de Cantabria. 39005-Santander, Spain \\
              \email{gila@unican.es}           
           \and
           Javier Segura \at
           Departamento de Matem\'aticas, Estad\'{\i}stica y Computaci\'on.
Facultad de Ciencias. Universidad de Cantabria. 39005-Santander, Spain.\\
\email{segurajj@unican.es} 
\and
Nico M.~Temme \at
IAA, 1825 BD 25, Alkmaar, The Netherlands. Former address: Centrum Wiskunde \& Informatica (CWI), 
        Science Park 123, 1098 XG Amsterdam,  The Netherlands.\\
\email{nicot@cwi.nl} 
}

\date{Received: date / Accepted: date}

\maketitle

\begin{abstract}

Methods for the
computation of classical Gaussian quadrature rules are described which are effective both for small and large degree. These methods are reliable because 
the iterative computation of the nodes has guaranteed convergence, and
they are fast due to their fourth-order convergence 
and its asymptotic exactness for an appropriate selection of the variables. 
For Gauss--Hermite and Gauss--Laguerre quadratures, local Taylor series can
be used for computing efficiently the orthogonal polynomials involved, with exact initial values for
the Hermite case and first values computed with a continued fraction for the Laguerre case. 
The resulting algorithms have almost unrestricted validity with respect to the parameters. 
Full relative precision is reached for the Hermite nodes, without any accuracy loss and for any degree,
and a mild accuracy loss occurs for the Hermite and Laguerre weights as well as for the Laguerre nodes. 
These fast methods are exclusively based on convergent processes, which, together with 
the high order of convergence
of the underlying iterative method, makes them particularly useful for high accuracy computations. We show examples of very high accuracy computations (of up to $1000$ digits of accuracy).
\end{abstract}

\section{Introduction}

Given a definite integral $I(f)=\int_{a}^{b}f(x)w(x)dx$, with $w(x)$ a weight function 
in an interval $[a,b]$, the $n$-point quadrature rule 
\begin{equation}
Q_n (f)=\displaystyle\sum_{i=1}^{n}w_i f(x_i)
\end{equation}
is said to be a Gaussian quadrature rule if it has the maximum possible degree of exactness, that is, if $I(f)=Q_n (f)$ for
$f$ any polynomial of degree not larger than the maximum possible degree, which is $2n-1$. 

As it is well known, the nodes $x_i$, $i=1,\ldots,n$ of the Gaussian quadrature rule are the roots of the (for instance monic) 
orthogonal polynomial satisfying 
\begin{equation}
\label{defpol}
\int_a^b x^{i} p_n (x) w(x) dx=0,\quad i=0,\ldots,n-1 .
\end{equation}
Among the Gauss quadrature rules, the most popular are those for which the associated orthogonal polynomials are solutions
of a linear second-order homogeneous ODE. These are the cases corresponding to 
classical orthogonal polynomials, namely:
Gauss--Hermite ($w(x)=e^{-x^2}$; $a=-\infty$, $b=+\infty$), 
Gauss--Laguerre ($w(x)=x^{\alpha} e^{-x}$, $\alpha>-1$; $a=0$, $b=+\infty$) and Gauss--Jacobi
($w(x)=(1-x)^{\alpha}(1+x)^{\beta}$, $\alpha,\beta>-1$; $a=-1$, $b=1$). The respective orthogonal
polynomials are denoted as $H_n (x)$ (Hermite polynomials), $L_n^{(\alpha)}(x)$ (Laguerre polynomials)
and $P_n^{(\alpha,\beta)}(x)$ (Jacobi polynomials).
All Gaussian quadratures with orthogonal polynomials satisfying a linear second-order homogeneous ODE are trivially related
to one of these three classical rules. In this paper we concentrate on Gauss--Hermite and Gauss--Laguerre
quadrature rules; Gauss--Jacobi quadrature will be described in a subsequent paper.

For the classical quadratures, the coefficients $a_n$, $b_n$ and $c_n$ of the three-term recurrence relation satisfied 
by the orthogonal polynomials, $P_{n+1}(x)=(a_n x + b_n) P_n (x) + c_n P_{n-1}(x)$, 
are available
in closed form 
and the nodes are
the eigenvalues of a tridiagonal matrix with entries in terms of the coefficients of the recurrence relation, 
while the weights can be computed from the eigenvectors  
 (see, for instance \cite[Section 5.3.2]{Gil:2007:NMF}). 
This procedure is generally  known as the Golub--Welsch algorithm 
\cite{Golub:1969:COG}, which was inspired by an observation made by Wilf 
\cite{Wilf:1978:MFT}. 
This is an interesting method for computing quadrature rules of low degree. However, as the number of nodes $n$
increases, the complexity scales as ${\cal O}\left(n^2\right)$ and the method slows down drastically.

An alternative to Golub--Welsch is the use of iterative methods, in which the central problem becomes the
computation of the nodes by some iterative root-finding method. This approach, which precedes Golub--Welsch in time (see, for instance, \cite{Davis:1956:AAW,Lowan:1942:TOT}), 
has recently received renewed attention, particularly for the computation of
high degree quadrature rules \cite{Glaser:2007:AFA,Hale:2013:FAA,Townsend:2016:FCO}. 
Some of the aforementioned iterative methods use asymptotic approximations (as the degree is large) 
for the nodes which are iteratively refined, with orthogonal polynomials also computed by
means of asymptotic expansions. This is the approach considered in \cite{Bogaert:2012:COL,JOH:2018:FAR}
for Gauss--Legendre quadrature, in \cite{Hale:2013:FAA} for Gauss--Jacobi quadrature and in 
\cite{Townsend:2016:FCO} for Gauss--Hermite; see also 
\cite{Swarztrauber:2002:OCT,Yakimiw:1996:ACW}. As an alternative numerical approach, 
we mention the recent work by Bremer \cite{Bremer:2017:OTC} on the computation of zeros 
of solutions of second-order ODEs via the computation of phase functions, which 
appears to be competitive for very large degrees ($10^5$ or larger).

As recently observed in \cite{Bogaert:2014:IFC}, only with asymptotic approximations it
is also possible to compute the nodes and weights of Gauss--Legendre quadrature in a non-iterative fashion, 
leading to very fast methods of computation. Similarly, it has been shown
in \cite{Gil:2018:GHL} that for Gauss--Hermite and Gauss--Laguerre a similar approach is possible. 
The same can be said regarding Gauss--Jacobi quadrature, as shown in \cite{GIL:2018:AEO} 
(which completes the asymptotic analysis of classical Gaussian quadratures). Both in 
\cite{Gil:2018:GHL,GIL:2018:AEO}, the validity of the expansions is limited to moderate values
of the parameters $\alpha$ and $\beta$. For other types of asymptotic approximations based on 
the Riemann--Hilbert approach, see \cite{Dea:2016:CAI,Huy:2018:CAI}.

Therefore, we have three main families of methods: the Golub--Welsch method, which is  an interesting
approach for low degrees; iteration-free asymptotic methods, which are preferable for large degrees;
 and iterative methods, which may provide the bridge between the two previous methods (particularly when
asymptotic estimations for large degrees are not used).

In this paper, we continue with the study of Gauss--Hermite and Gauss--Laguerre quadratures initiated
in \cite{Gil:2018:GHL}, and we now consider purely iterative methods which are
 free of asymptotic approximations but which are asymptotically exact, in the sense that 
for large
degrees the number of iterations required per node tends to $1$.

With respect to the Golub--Welsch algorithm, our method is particularly advantageous as the degree becomes large,
as is also the case of the other aforementioned iterative methods. 
With respect to previous iterative methods, the present method has the crucial advantage of its higher rate 
of convergence, its reliability (convergence is proved) and its larger range of applicability (almost unrestricted). And with respect to the iterative methods based on asymptotics \cite{Hale:2013:FAA,Townsend:2016:FCO}, it has the additional
advantage that arbitrary precision is available, and for any value of the parameters 
(small or large degrees, and unrestricted $\alpha$ for Gauss--Laguerre).
This last advantage with respect to iterative-asymptotic methods also holds with respect to purely asymptotic methods as those in \cite{Gil:2018:GHL,GIL:2018:AEO}.

We expect that an optimal
algorithm for the computation of Gauss quadratures in fixed precision will involve both the 
asymptotics-free iterative methods and the iterative-free asymptotic methods, probably complemented with
the Golub--Welsch algorithm for small degree. The present paper is a necessary step in this direction.

\section{Reliable iterative computation of Gaussian quadratures}

In this section we describe the main general ingredients in the iterative computation of Gauss--Hermite
and Gauss--Laguerre quadrature rules; in later sections  we analyze the particular methods used for computing
the orthogonal polynomials involved as well as associated values (weights) both for the
Gauss--Hermite (Section~\ref{Gauss--Hermite}) and Gauss--Laguerre (Section~\ref{Gauss--Laguerre}) quadratures. 
In this section we first summarize briefly the main ingredients of the 
fourth-order fixed point method \cite{Segura:2010:RCO}, which will be our choice of iterative method for 
solving non-linear equations. Then, we study the possible Liouville transformations of the ODEs which result 
in different possible selections of globally convergent fixed point methods. 
After this, we consider the computation of the weights in terms of the derivatives of the Liouville-transformed 
functions, and we show that for a particular change of variables (which we call canonical), 
the method is asymptotically exact for the most significant nodes; we obtain well-conditioned expressions 
for the weights in terms of the canonical variable. Finally, we outline the method of computation of 
orthogonal polynomials (or related functions), which is later explained in more detail for the Hermite
and Laguerre cases (Sections  \ref{Gauss--Hermite} and  \ref{Gauss--Laguerre}).

\subsection{The iterative method}
\label{itera}

There has been an almost general consensus in using Newton's method as iterative method for computing
the nodes of Gaussian quadratures, and only in \cite{Yakimiw:1996:ACW}
a higher order variant is considered. 
Newton's method is a well-known generic method for solving non-linear equations. However, for the particular 
case of functions
which are solutions of second-order ODEs, better methods exist. 
In particular, the method introduced in \cite{Segura:2010:RCO} 
has essentially the same computational cost as Newton's method but it has three fundamental advantages: 
it doubles the order of convergence of Newton's method, it 
converges with certainty
and, as commented before, with the 
appropriate selection of variable, the method 
tends to be exact as the degree goes to infinity (it gives the exact root in one step).

This fixed point method is able to compute all the
zeros of any solution of a second-order ODE in normal form (without first derivative term) 
$y''(x)+A (x) y(x)=0$ provided that  $A(x)$ is continuous and the monotonicity properties 
of $A(x) $ in this interval are known in advance. No initial estimations of the zeros are needed, and
the method computes all the zeros with certainty in the direction of decreasing values of $A(x)$.

For the moment, we assume that a method for computing function
values for the classical orthogonal polynomials is available.

The equation being in normal form is not an important restriction, because any differential equation
\begin{equation}
\label{ODE2}
\label{normf}
w''(x)+b(x)w'(x)+a(x)w(x)=0,
\end{equation}
with $b(x)$ differentiable can be transformed into normal form with a change of function, 
a change of variables $z=z(x)$ or both (see next subsection).

This fixed point method can be understood as a consequence of the following Sturm theorem:

\begin{theorem}[Sturm comparison]
Let $y(x)$ and $v(x)$ be solutions of $y''(x)+A_{y}(x)y(x)=0$ and $v''(x)+A_{v}(x)v(x)=0$ respectively, with
 $A_v (x)>A_{y}(x)$. If $y(x^{(0)})v'(x^{(0)})-y'(x^{(0)})v(x^{(0)})=0$ 
and $x_y$ and $x_v$ are the real zeros of $y(x)$ and $v(x)$
closest to $x^{(0)}$ and larger (or smaller) than $x^{(0)}$, then $x_v<x_y$ (or $x_v>x_y$).
\end{theorem}

This theorem is easy to prove and has a simple geometrical interpretation in terms of the speed of oscillation of the
solutions, which is greater as the coefficient of the ODE becomes greater. See for instance \cite{Gil:2014:RSD}.

As a consequence of this theorem, a method for the computation of the zeros
of solutions of $y''(x)+A (x)y(x)=0$ emerges. If $A (x)$ is a decreasing (increasing)
function and $A (x)>0$, we compute the zeros with an increasing (decreasing) sequence (if $A(x)<0$ in an interval 
the solutions have one zero at most in this interval). Given a value $x^{(0)}$, 
the zero of $y(x)$ closest to $x^{(0)}$ and larger (smaller) than $x^{(0)}$ can be computed with 
certainty using the following scheme.

\begin{algorithm}[Zeros of $y''(x)+A (x)y(x)=0$, $A (x)>0$ monotonic].
\label{Algo1}

Let $x^{(0)}<\alpha$ with $y(\alpha )=0$ and such that there is no zero of $y(x)$ 
between $x^{(0)}$ and $\alpha$, and assume that $A (x)$ is decreasing (increasing). 

Starting from $x^{(0)}$, 
compute $x^{(n+1)}$ from $x^{(n)}$ as follows: 
find a non-trivial solution of the equation $v''(x)+A (x^{(n)})v(x)=0$ such that
$y(x^{(n)})v'(x^{(n)})-y'(x^{(n)})v(x^{(n)})=0$. Take as $x^{(n+1)}$ 
the zero of $v(x)$ closest to $x^{(n)}$ and larger (smaller) than 
$x^{(n)}$. Then, the sequence $\{x^{(n)}\}$ converges 
monotonically to $\alpha$.

\end{algorithm}

  Observe that solving the differential equation $v''(x)+A (x^{(n)})v(x)=0$ 
of the previous theorem is trivial because the coefficient is constant.

 The algorithm can be applied successively to generate a sequence of zeros as follows.
\begin{algorithm}[Computing a sequence of zeros, $A (x)$ monotonic]
\label{Algo2}

Let $\alpha_1, \alpha_2$ be consecutive zeros of $y(x)$, with $\alpha_1 <\alpha_2$. 

If $A (x)$ is decreasing and $\alpha_1$ is known, the zero $\alpha_2$ can be computed using
Algorithm \ref{Algo1} with starting value $x^{(0)}=\alpha_1$ 
(the first iteration being $x^{(1)}=\alpha_1+\pi/\sqrt{A (\alpha_1)}$).

 If $A (x)$ is increasing and $\alpha_2$ is known, the zero $\alpha_1$ can be computed using
Algorithm \ref{Algo1} with starting value $x^{(0)}=\alpha_2$ 
(the first iteration being $x^{(1)}=\alpha_2-\pi/\sqrt{A (\alpha_2)}$).
\end{algorithm}

As commented before the sequences generated are increasing (decreasing) if $A (x)$ is decreasing (increasing).

The iteration of Algorithm \ref{Algo1} can be explicitly written as follows:

\begin{equation}
\label{def1}
T_j (x)=x-\Frac{1}{\sqrt{A(x)}}\arctan_{j}\left(\sqrt{A(x)}\,h(x)\right),
\end{equation}
with $h(x)=y(x)/y'(x)$, $j=\mbox{sign}(A'(x))$ and
\begin{equation}
\label{def2}
\arctan_{j}(\zeta)=\left\{
\begin{array}{l}
\arctan(\zeta)
\mbox{ if }j\zeta > 0,\\ 
\arctan(\zeta)+j\pi \mbox{ if }j\zeta\le 0,\\ 
j\pi/2 \mbox{ if } \zeta=\pm \infty.
\end{array}
\right.
\end{equation}
Observe that the only fixed points of $T_j(x)$ are the zeros of $y(x)$.

The algorithms need some a priori analysis: the monotonicity properties of the coefficient $A (x)$ 
must be known in advance, because
the method has to be applied separately in those subintervals where $A (x)$ is monotonic.
This analysis has been completed for hypergeometric functions \cite{Deano:2004:NIF}, and we will use this
information in our algorithms.

Initial estimations are not needed, but as we will see the use of 
some simple bounds for the extreme zeros  \cite{Dimitrov:2010:SBF} in order to refine the stopping criterion is convenient.

From the construction of the method, we observe that it is exact (gives the exact roots in one iteration) if $A (x)$ is constant.
Therefore, if in some limit the coefficient of the ODE tends to a constant value, then the method is 
asymptotically exact in that limit. 

\subsection{Liouville transformations of the differential equations and computation of the nodes}
\label{Lio}

The classical orthogonal polynomials $H_n(x)$, $L_n^{(\alpha )}(x)$ and $P_n^{(\alpha,\beta)}(x)$ satisfy 
second-order ODEs
(\ref{ODE2})
with $a(x)$ and $b(x)$ simple rational coefficients. The iterative method described in Section
\ref{itera} requires that the ODE is in normal form (\ref{normf}); in addition, 
the method requires that the monotonicity properties of the coefficient of the ODE are known in advance.
The ODEs for orthogonal polynomials can be transformed into their normal forms by Liouville transformations in which the 
changes of variables can be selected conveniently in order to simplify the analysis of the coefficient. The 
necessary analysis was performed in \cite{Deano:2004:NIF,Deano:2007:GSI}.

 Given a function $w(x)$ which is a solution of Eq.~(\ref{ODE2}) and a change of the independent variable
$z=z(x)$, then the function $y(z)$, with $y(z(x))$ given by

\begin{equation}
y(z(x))=\sqrt{z'(x)}\exp\left(\frac12\displaystyle{\int}^{x} b(x)\right) w(x) ,
\label{funcL}
\end{equation}
\noindent satisfies the equation in normal form 
 \begin{equation}
 \ddot{y}(z)+A (z) y(z)=0\,,
 \label{normal}
 \end{equation}
 \noindent where the dots represent differentiation with respect to $z$ and 

 \begin{equation} 
  A (z) =\dot{x}^2  \tilde{A}(x(z))+\frac12 \{x,z\},\quad \tilde{A}(x)=a- b'/2-b^2/4,
 \end{equation}
 \noindent where 
 $\{x,z\}$ is the Schwarzian derivative of $x(z)$
 with respect to $z$ \cite[p. 191]{Olver:1997:AAS}.
  As a function of the original variable $x$ this can be written
 \begin{equation}
 \begin{array}{ll}
 A(z(x)) & = \Frac{1}{z'(x)^2}(\tilde{A}(x)
 -\frac12\{z,x\})\\
&  
= \Frac{1}{d(x)^2}\left(a(x)-\Frac{b'(x)}{2}-\Frac{b(x)^2}{4} + \Frac{3d'(x)^2}
{4d(x)^2}-\Frac{d''(x)}{2 d(x)} \right),
 \end{array}
 \label{master}
 \end{equation}
\noindent where $\{ z ,x\} $ is the Schwarzian derivative of $z(x)$ 
with respect to $x$ and $d(x)=z'(x)$. 
 
In \cite{Deano:2004:NIF,Deano:2007:GSI} a systematic study  of the Liouville transformations that 
lead to second-order equations with simple enough coefficients $A(x)$ was performed for the confluent and
Gauss hypergeometric equations (and therefore, in particular for classical orthogonal polynomials). We briefly 
describe the cases for Hermite and Laguerre concentrating on the changes
of variables most useful for our purpose. Gauss--Jacobi rules will be described in detail in a future publication, 
and we will advance at the end of this paper some ideas about those rules.

\subsubsection{Hermite polynomials}
\label{AlgHe}

The function $w(x)=H_n(x)$ satisfies the ODE
$$
w''(x)-2x w'(x)+2n w (x)=0.
$$
We transform to normal form without changing the variable $x$, and write $y(x)=e^{-x^2/2}H_n (x)$. This function satisfies
\begin{equation}
\label{normH}
y''(x)+A(x) y(x)=0,\,A(x)=2n+1-x^2.
\end{equation}
The coefficient $A(x)$ is very simple and no change of variables is needed. In addition, as $n$ becomes large
the coefficient becomes approximately constant for small $x$; this means that the fixed point method will improve 
its convergence speed as $n$ becomes large, particularly for the small zeros, which, as we will, see are
the most significant nodes (those with the largest weights). 

In this case, because of the symmetry of the zeros, we only need to consider the positive zeros. The fixed point method proceeds starting from $x=0$ and computing zeros in the direction of increasing $x$ (decreasing $A(x)$), 
which is the direction of decreasing weights.

The methods do not need sharp estimations for the roots. The method terminates when $\lfloor n/2 \rfloor$ positive roots
are obtained.

\subsubsection{Laguerre polynomials}
\label{Lagcha}

The function $w(x)=L_n^{(\alpha)}(x)$ satisfies
$$
w'' (x)+\left(\Frac{\alpha+1}{x}-1\right)w'(x)+\Frac{n}{x}w(x)=0.
$$
Without a change of the variable $x$, we transform to normal form and we obtain
\begin{equation}
\label{norlx}
y(x)=x^{(\alpha+1)/2}e^{-x/2}w(x),
\end{equation}
which satisfies
\begin{equation}
\label{Alx}
y''(x)+A(x)y(x)=0,\quad A(x)=\Frac{1}{4}\left(-1+\Frac{2L}{x}+\Frac{1-\alpha^2}{x^2}\right),
\end{equation}
where $L=2n+\alpha+1$.
The coefficient is simple as also the monotonicity properties are
(decreasing if $|\alpha|<1$ and with a maximum at 
$(\alpha^2-1)/L$ if $\alpha >1$).

As shown in \cite{Deano:2004:NIF}, the changes $z(x)=\frac{1}{m}x^m$ for $m\neq 0$ and $z(x)=\log(x)$ give Liouville transformations
which also lead to ODEs in normal form with at most one 
extremum of the resulting coefficient (except for some cases when 
$m\in (0,1/2)$). The resulting differential equation
$$
\ddot{y}(z)+A(z)y(z)=0
$$
(where dots mean derivatives with respect to $z$) is such that
$$
A(z(x))=\frac14 x^{-2m} (-x^2+2Lx+m^2-\alpha^2), 
$$
where the case $z(x)=\log x$ corresponds to $m=0$. We have an infinite number of possible changes available, but an interesting
selection is $m=1/2$ because, as happened with the Hermite case, we have a constant term which grows 
with $n$, which is interesting from the point of view of the asymptotic exactness of the method 
as $n\rightarrow +\infty$. 
Then we take 
$z(x)=\frac{1}{m}x^m$ for $m=1/2$ or equivalently $z(x)=\sqrt{x}$, and we have that
\begin{equation}
\label{yzl}
y(z)=z^{\alpha+1/2} e^{-z^2/2} L_n^{(\alpha)}(z^2)
\end{equation}
satisfies
\begin{equation}
\label{normL}
 \ddot{y}(z)+A(z)y(z)=0,\quad A(z(x))=-x+2L+\frac{\frac14- \alpha^2}{x},
\end{equation}
and $A(z(x))$ is decreasing for positive $x$ if $|\alpha|\le 1/2$ and has a maximum at $x_e=\sqrt{\alpha^2 -1/4}$ if $|\alpha|>1/2$.
The fixed point method can therefore be applied to the function (\ref{yzl}) with ease.

No initial estimations for the roots are required. However, it is convenient to use bounds for the extreme zeros in order
to stop the method. Then, if $|\alpha|\le 1/2$, we can start the process from $z$ equal to square root of 
the lower bound for the zeros
(see \cite[Eq.~(1.2)]{Dimitrov:2010:SBF}). If $|\alpha|>1/2$ we start from the maximum of $A(z)$, which is $z_e=(\alpha^2-1/4)^{1/4}$,
and compute zeros in increasing order until the upper bound is surpassed 
(because all the values of $z$ generated constitute a monotonic sequence, this is a safe stopping rule); 
after this, we start again from $z_e$ and compute zeros in the direction of 
decreasing $z$ until a total of $n$ zeros has been computed.

\subsection{Computation of the weights}

As before, in this section we assume that an algorithm for the computation of the orthogonal polynomials is available
(we discuss in Sections \ref{Gauss--Hermite} and \ref{Gauss--Laguerre} how to compute them). 
We now describe the computation of the weights assuming that the nodes have already been computed.

\subsubsection{Gauss--Hermite weights}
\label{ghw}

As it is well known, in terms of the first derivative, the Gauss--Hermite weights can be written as
\begin{equation}
w_i = \Frac{\sqrt{\pi}2^{n+1}n!}{[H_n ^{\prime}(x_i)]^2}.
\end{equation}
Considering now the solution of (\ref{normH}), $y(x)=e^{-x^2/2}H_n(x)$. In terms of this function, the weights become
\begin{equation}
w_i=\Frac{\sqrt{\pi}2^{n+1}n!}{[y^{\prime}(x_i)]^2}e^{-x_{i}^2}\equiv \omega_i e^{-x_{i}^2},
\end{equation}
and we say that $\omega_i$ are the scaled weights.

Observe that, because the coefficient of the ODE (\ref{normH}) is essentially constant when $n$
is large, then $|y^{\prime}(x)|$ should be approximately constant, and the main dependence
on the nodes is in the exponential factor $e^{-x_{i}^2}$. This is confirmed using asymptotics
for $n\rightarrow +\infty$, which gives
\begin{equation}
\label{asweh}
w_i \sim \Frac{\pi}{\sqrt{2n}}e^{-x_i^2},
\end{equation}
where the estimation works better for the small zeros. With this relation, we observe that
the weights decrease exponentially as we move away from $x=0$ which, as explained in the previous
section, is the starting point for the fixed point method; this
method will compute nodes in the direction of decreasing weights, starting from the most significant
nodes. And for these first nodes the method is more rapidly convergent as $n$ becomes larger (asymptotic 
exactness).

\subsubsection{Gauss--Laguerre weights} 
\label{glw}

In terms of the first derivative, the Gauss--Laguerre weights are
\begin{equation}
w_i = \Frac{\Gamma (n+\alpha+1)}{n! x_i \left[L_n ^{(\alpha) \prime}(x_i)\right]^2}=\Frac{4\Gamma (n+\alpha+1)}{n! 
\left[\Frac{d}{dz}L_n ^{(\alpha)}(z_i ^2)\right]^2},
\end{equation}
where, as in the previous section, $x=z^2$.

In terms of (\ref{yzl}), solution of (\ref{normL}),
\begin{equation}
\label{wasym2}
w_i = \Frac{4\Gamma (n+\alpha+1)}{n! \left[\dot{y}(z_i)\right]^2}x_i^{\alpha+\frac12}e^{-x_i}
\equiv \omega_i x_i^{\alpha+\frac12}e^{-x_i},
\end{equation}
where the $\omega_i$ are the scaled weights.

As for the Hermite case, the coefficient of the ODE (\ref{normL}) is essentially constant when $n$
is large, particularly around its maximum when $|\alpha|>1/2$ (at $x_e=\sqrt{\alpha^2-1/4}$), and 
then $|\dot{y}(z_i)|$ should be approximately constant; the main dependence
on the nodes is in the exponential factor $x_i^{\alpha+1/2}e^{-x_i}$. Again, this is confirmed 
considering asymptotic estimates as $n\rightarrow \infty$:
\begin{equation}
\label{wasym}
w_i  \sim \pi\Frac{\Gamma (n+\alpha +1)}{n!}n^{-\alpha-\frac12}x_i ^{\alpha+\frac12} e^{-x_i}\sim \pi n^{-1/2} x_i ^{\alpha +\frac12} e^{-x_i},
\end{equation}
where the estimation works better for the small zeros.

Observe that the function $f(x)=x^{\alpha + 1/2}e^{-x}$ has its maximum at $x_M=\alpha+1/2$ when
$\alpha>-1/2$ and that this will be close to the starting point for the fixed point method, which
is $x_e=\sqrt{\alpha^2-1/4}$ when $|\alpha|>1/2$. Then, as happened for the Hermite case, in the 
canonical variable (which is $z=\sqrt{x}$ for Laguerre) the fixed point
method will compute nodes in the direction of decreasing weights, starting from the most significant
nodes (also for $|\alpha|<1/2$, because, in this case,  the first computed node is the smallest). 
And for these nodes the method is more rapidly convergent as $n$ becomes larger (asymptotic 
exactness), because they are close to the maximum of $A(z(x))$ when $|\alpha|>1/2$.

The fact that the method computes first the most significant nodes (and faster as the
degree increases due to asymptotic exactness) and successively the rest of nodes in decreasing order, is also 
interesting if subsampling is to be considered, that is, if only the nodes with weights larger than 
a given threshold are of interest.

\subsubsection{Scaled weights: condition and range of computation}
\label{scaw}

The scaled weight for the Hermite case can be written as $\omega_i =\omega(x_i)$, with 
$\omega (x) = k_n |y'(x)|^{-2}$, where $y(x)$ is a solution of the second-order ODE (\ref{normH}) 
and $k_n$ only depends on $n$. Similarly, for Laguerre $\omega_i =\omega(z_i)$, with 
$\omega (z) = k_{n,\alpha} |\dot{y}(z)|^{-2}$, 
where $y(z)$ is a solution of the second-order ODE (\ref{normL}) 
and $k_{n,\alpha}$ only depends on $n$ and~$\alpha$.

We note that the scaled weights are well conditioned as a function of the nodes in the canonical variable.
This is so because, considering for instance the Hermite case, $\omega' (x)= -2 c_n y'(x)^{-3} y''(x)$,
but $y''(x_i)=0$ because $y(x_i)=0$ and $y(x)$ satisfies an equation in normal form; therefore 
$\omega' (x_i)=0$ and at first
order the scaled weights do not depend on the values of the nodes. The same is true for the Laguerre case
in terms of the $z$ variable. The main source of errors in the computation of the weights is in the 
elementary function which has been factored out for the scaled weights.

Both for the Hermite and Laguerre cases, the computation of the scaled weights is free 
of overflow/underflow problems, both as a function of the nodes and the parameters. The main dependence on the nodes is
factored out in an elementary function, while the dependence of the scaled weights on the degree
goes as $n^{-1/2}$ for $n$ large, and the dependence on $\alpha$ will be of no concern, as we 
explain next. 

It is in fact possible to compute the scaled weights without computing the constants $k_n$ and
$k_{n,\alpha}$. Considering the Hermite case (the same idea works for Laguerre), 
the idea is to solve the ODE (\ref{normH}) with some arbitrarily
chosen normalization for the solutions; then compute $\tilde{\omega}_i=1/y'(x_i)^2$, which are proportional
to the scaled weights. Finally, the constant of proportionality can be fixed by using the fact that the sum of the
(unscaled) weights is $\sqrt{\pi}$ for Gauss--Hermite (and $\Gamma (\alpha +1)$ for Laguerre). 
Proceeding in this way, we eliminate possible overflows/underflows with respect to the degree $n$, and also
with respect to the parameter $\alpha$ for the Laguerre case.
This leads to practically unrestricted 
algorithms for scaled weights, while for the original weights the possible underflows are controlled by
an elementary factor. 

The only issue which remains to be discussed is how the orthogonal polynomials are computed. The approach
varies depending on the type of quadrature. For the Hermite case we can solve 
the problem just by using Taylor series, while for the Laguerre case Taylor series should be
supplemented with a continued fraction evaluation.

\section{Computing the Gauss--Hermite quadrature}
\label{Gauss--Hermite}

The Gauss--Hermite quadrature has the special property that the differential
equation does not have finite singularities; the same is not true for the other
classical Gauss quadratures. This enables the possibility
of computing the polynomials by local Taylor series, similarly as done 
in \cite{Glaser:2007:AFA}. For other cases, and in particular for Gauss--Laguerre, 
local Taylor series are also possible away
from the singularities, but the application of series is necessarily
more limited and must be complemented with other methods.

The algorithm consists in the computation of the nodes with the fixed point method described
in Section~\ref{itera} and as described in Section~\ref{AlgHe}, with weights
computed following Sections \ref{ghw} and \ref{scaw}. The function 
$y(x)=\lambda e^{-x^2/2}H_n (x)$, with $\lambda$
a constant which is introduced for later convenience, and
its derivative are computed in parallel with the application of the fixed point method. We 
describe the method step by step.

The nodes are symmetric around the origin, and $x=0$ is a zero for odd degree, and, as
described in Section~\ref{AlgHe}, the fixed point method starts from $x=0$ and computes 
zeros in the direction of increasing $x$ (decreasing $A(x)$), 
which is the direction of decreasing weights.

We start at $x^{(0)}=0$, and the first step of the algorithm is

$$
x^{(1)}=T_{-1}(x^{(0)})=\left\{
\begin{array}{lll}
\Frac{\pi}{\sqrt{2n+1}}&,& n \mbox{ odd},\\
\Frac{\pi}{2\sqrt{2n+1}}&,& n \mbox{ even}.
\end{array}
\right.
$$

Observe that $h(0^+)=y (0^+)/y^{\prime}(0^+)=0^+$ is $n$ is odd (and $x=0$ is a node) and 
$h(0^+)=y (0^+)/y^{\prime}(0^+)=+\infty$ if $n$ is even (and $x=0$ is not a node).
Notice that this value $x^{(1)}$ is a lower bound for the first positive zero.

For computing the second step, $x^{(2)}=T_{-1}(x^{(1)})$, we need to compute $y(x^{(1)})$ and
 $y'(x^{(1)})$, and for this we use Taylor series centered at $x^{(0)}=0$ 
(using the known values $y (0)$ and $y^{\prime}(0)$). A way to choose these values,
 without needing to compute $H_n(0)$ and $H^{\prime}_n (0)$ (which involves computing
some factorials) is selecting initially an arbitrary normalization and then rescaling at the end; 
this is the approach we consider in our algorithms. 
We will take $y (0)=1$, $y^{\prime}(0)=0$ if $n$ is even, and $y (0)=0$, 
$y^{\prime}(0)=1$ if $n$ is odd.

The truncated local Taylor series are
\begin{equation}
\label{LTS}
y (x+h)=\displaystyle\sum_{i=0}^{N}\Frac{y^{(i)}(x)}{i!}h^i,\quad
y^{\prime} (x+h)=\displaystyle\sum_{i=0}^{N}\Frac{y^{(i+1)}(x)}{i!}h^i,
\end{equation}
where, for the first step, $x=x^{(0)}$ and $h=x^{(1)}-x^{(0)}$. 
The derivatives $y^{(i)}(x)$ can be computed by differentiating the ODE satisfied by 
$y (x)=A e^{-x^2 /2}H_n(x)$. We have
\begin{equation}
\label{reH}
y^{(k+2)}(x)+(2n+1-x^2)y^{(k)}(x)-2kxy^{(k-1)}(x)-k(k-1)y^{(k-2)}(x)=0,
\end{equation}
and we feed this recurrence relation with the known values $y(x)$ and $y'(x)$. In the truncated series, the value of $N$ does not need to be fixed a priori, and we can sum the series until the last term gives a relative contribution
smaller than the relative accuracy goal.

The algorithm proceeds in the same way in each iteration. 
After it computes a new iteration, $x^{(i)}=T_{-1}(x^{(i-1)})$, the values
of $y(x^{(i)})$ and $y'(x^{(i)})$, needed to compute the next iteration  
$x^{(i+1)}=T_{-1}(x^{(i)})$, are evaluated by using Taylor series centered at $x_{i-1}$ with step
$h=x^{(i)}-x^{(i-1)}$. The process is repeated until an accurate approximation to the first
positive node $\alpha_1$ is obtained. Then, we start a new iterative process for the next zero
$\alpha_2>\alpha_1$ with $x^{(0)}=\alpha_1$, and 
$x^{(1)}=T_{-1}(\alpha_1)=\alpha_1+\pi / \sqrt{A(\alpha_1)}$, and iterate until convergence
to $\alpha_2$ is reached; and so on. The process can be stopped after $\lfloor n/2 \rfloor$ positive 
nodes have been computed. In this process, the approximate values of $y'(\alpha_i)$ are also stored,
and they will be used to compute the nodes, as we are going to explain. Before this, we first 
comment on the stability of the recursive process to compute derivatives~(Eq.~(\ref{reH})).

Eq.~(\ref{reH}) is a difference equation of fourth-order (a five-term recurrence relation), and,  therefore, the linear space of solutions 
of this recurrence relation has dimension $4$. For the forward computation of derivatives, it is 
essential that the derivatives of $y(x)$ are not recessive as $n\rightarrow \infty$, which means that
there are no other solutions of the recurrence relation (\ref{reH}), say $g_n$, such that 
$\displaystyle\lim_{n\rightarrow \infty} y_n /g_n =0$; if that were the case the computation would be unstable.
The Perron-Kreuser theorem \cite{Kreuser:1914:UDV} is a simple tool to analyze the conditioning of linear recurrence
relations (see \cite{Cash:1980:ANO} for a more recent account of this result). For the case of (\ref{reH}) this theorem is not
conclusive and it gives the information that all solutions of this difference equation satisfy 
$\limsup_{k\rightarrow +\infty}\left(|y^{(k)}|/(k!)^{2/3}\right)^{1/k}=1$. This, on one hand, indicates that the radius
of convergence of local Taylor series (Eq.~(\ref{LTS}) with $N=\infty$) is infinity, which could be expected given that the ODE has no finite singularities. On
the other hand, the fact that all the solutions have this behaviour means that there are no solutions 
of (\ref{reH}) 
that are exponentially larger than other ones as $n\rightarrow \infty$. This suggests, although it does not imply stability,  that the computation of derivatives can be stable, as numerical experiments indeed confirm.

The initial values we have considered, ($y (0)=1$, $y^{\prime}(0)=0$ if $n$ is even, and $y (0)=0$, 
$y^{\prime}(0)=1$ if $n$ is odd), together with the fact that in each step 
we are integrating the differential equation (\ref{normH}), means that we are computing values
of $y(x)=\lambda e^{-x^2/2} H_n (x)$ and its derivative, with $\lambda$ an unknown constant that
should be evaluated in order to compute the weights; we describe next how to fix this normalization. 

From the approximate values at the positive 
nodes $\alpha_i$, we can compute
the quantities $\bar{\omega}_i=|y'(\alpha_i)|^{-2}$, $i=1,\ldots \lfloor n/2\rfloor$ 
and these quantities will be proportional to the corresponding scaled 
weights $\omega_i$.\footnote{Notice that we have changed the notation for the nodes with respect to the previous section and
 the index runs differently 
since we are considering the positive nodes.} Observe that, if $n$ is odd, with our normalization
the value corresponding to the node $\alpha_0=0$ is $\bar{\omega}_0=1$.
Then we have $\tilde{\omega}_i=C \omega_i$, where
$C$ is constant that can be fixed by computing one of the moments.
For instance, we have that, 
denoting as before by $\alpha_1<\alpha_2<\ldots$ the positive nodes 
and $w_1,w_2,\ldots$ their
corresponding weights
\begin{equation}
\label{normaw}
\mu_1=\displaystyle\int_{-\infty}^{\infty} x^2 e^{-x^2}dx=\Frac{1}{2}\sqrt{\pi}=2\displaystyle\sum_{j=1}^{\lfloor n/2 \rfloor}
w_j \alpha_j^2 .
\end{equation}
We use this normalization to fix the correct normalization for the 
weights.\footnote{Of course, we could use other moments, as for instance 
$\displaystyle\int_{-\infty}^{\infty} e^{-x^2}dx=\sqrt{\pi}=
\delta_{1,k}w_0 f(0)+2\displaystyle\sum_{j=1}^{\lfloor n/2 \rfloor}
w_j f(\alpha_j)$, $k=n-2\lfloor n/2\rfloor$.}

From the values of $\bar{\omega}_i$ obtained, we compute $\bar{w}_i=
\bar{\omega}_i e^{-\alpha_i^2}$, and we have that
\begin{equation}
\label{normaw2}
\bar{\mu}_1=2\displaystyle\sum_{j=1}^{\lfloor n/2 \rfloor}
\bar{w}_j \alpha_j^2 =C\mu_1=C\sqrt{\pi},
\end{equation}
from where we obtain $C$, the scaled weights $\omega_i=\bar{\omega}_i /C$ and
the weights $w_i=\bar{w}_i /C$. Notice that few terms will be needed 
in the sum (\ref{normaw2}), even when $n$ is large, because of the exponential
decay of the weights. With this procedure to compute the weights, overflow/underflow
problems are completely eliminated for the scaled weights and we have explicit 
control of the main dominant exponential factor for the unscaled weights $w_i$.

\subsection{Numerical results}
\label{numres}

The resulting algorithm is short and simple and very efficient. The only ingredients are the application of the 
fixed point iteration, the use of truncated Taylor series and the normalization (\ref{normaw}). No accuracy degradation takes places for the nodes and full accuracy is reached (but some mild error degradation does take place for the 
weights, as we later discuss). This is in 
contrast to what is described in \cite{Glaser:2007:AFA} and 
\cite{Townsend:2016:FCO}, and is certainly a notable property of the method. We note that the nodes are computed
in increasing order, and that therefore this is favourable for the stability in the computation of the nodes; 
that the algorithm in \cite{Glaser:2007:AFA}, which computes nodes in the same order, 
accumulates some error in the nodes could be a consequence of the fact that additional techniques are needed which are absent in our algorithm (like first estimations of the nodes using a Runge--Kutta method). Exploring the differences
in error propagation should involve a detailed comparison between the methods and implementations; in this 
sense it is worth noticing that very subtle differences may influence error propagation, as we later discuss 
(see Section~\ref{subtle}).

\begin{figure}
\vspace*{0.8cm}
\centerline{\includegraphics[height=6cm,width=12cm]{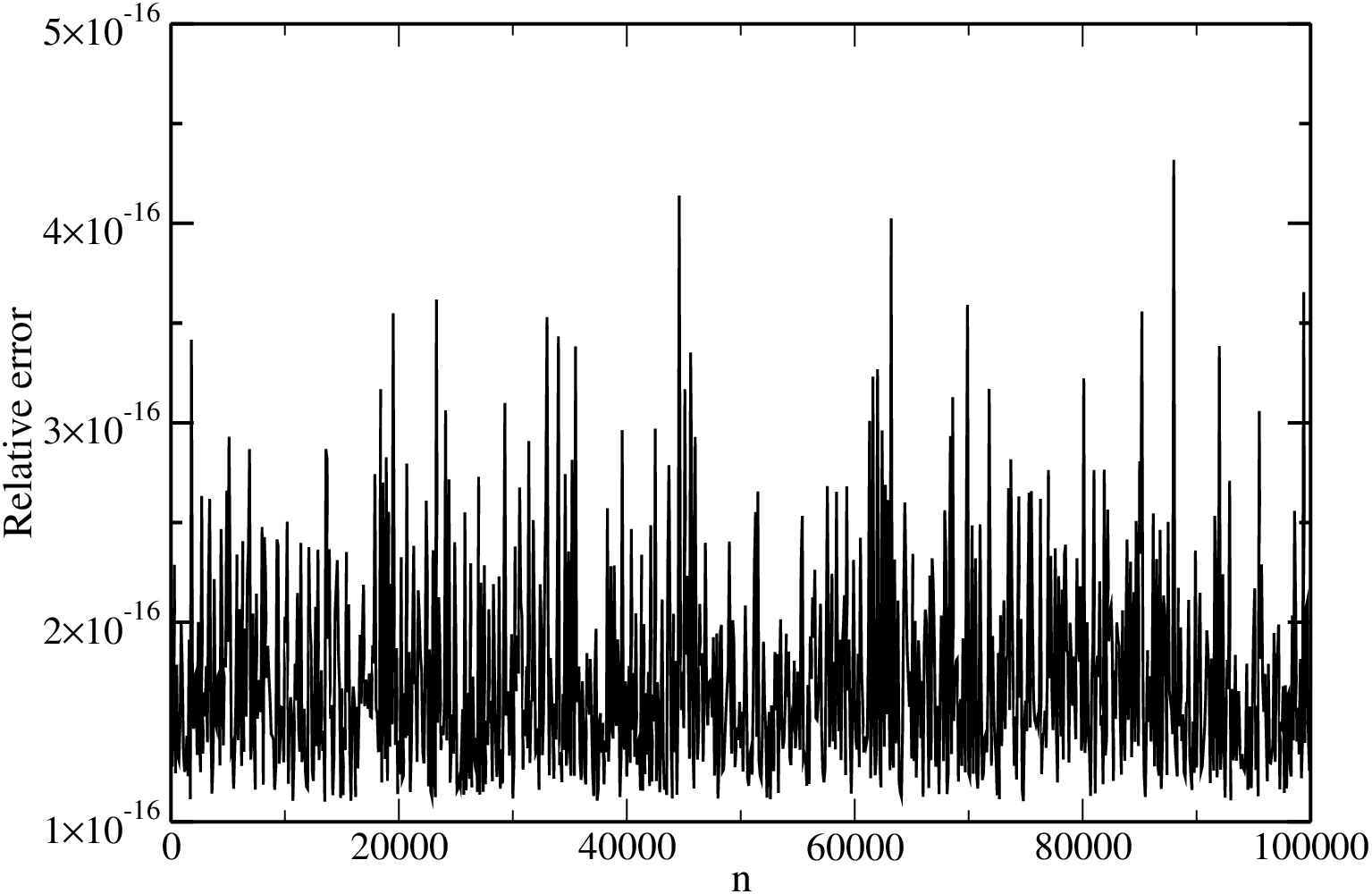}}
\caption{Relative errors in the computation of the nodes for $n$-point Gauss--Hermite
quadrature. For each $n$, the value 
$\displaystyle\max_{i=1,\ldots n}|1-x_{i}^{(d)}/x_i^{(q)}|$
is represented, where $x_i^{(d)}$ are the nodes computed in double precision and 
$x_i^{(q)}$ are the nodes in quadruple precision (the trivial node $x=0$ for
odd degree is excluded).
}
\label{fig1}
\end{figure}

Figure~\ref{fig1} shows the maximum relative errors in the 
computation of the nodes for orders $n$ smaller than $10^5$; 
the errors are obtained by comparing the nodes 
obtained with our algorithm in double precision (coded in Fortran) with a 
quadruple precision version of the same algorithm. The figure shows the typical 
noise pattern consistent with double precision accuracy and a detailed inspection
shows that all digits are correct except, in some cases, the last digit and by a
small amount. This is a surprising result, and it is in part explained by the fact that the zeros are computed
in increasing order, but it is not the only reason (in \cite{Glaser:2007:AFA} 
the zeros are also computed in increasing order, but some error
degradation happens); for additional details on this notable behaviour, see Section~\ref{subtle}.

\begin{figure}
\vspace*{0.8cm}
\centerline{\includegraphics[height=6cm,width=12cm]{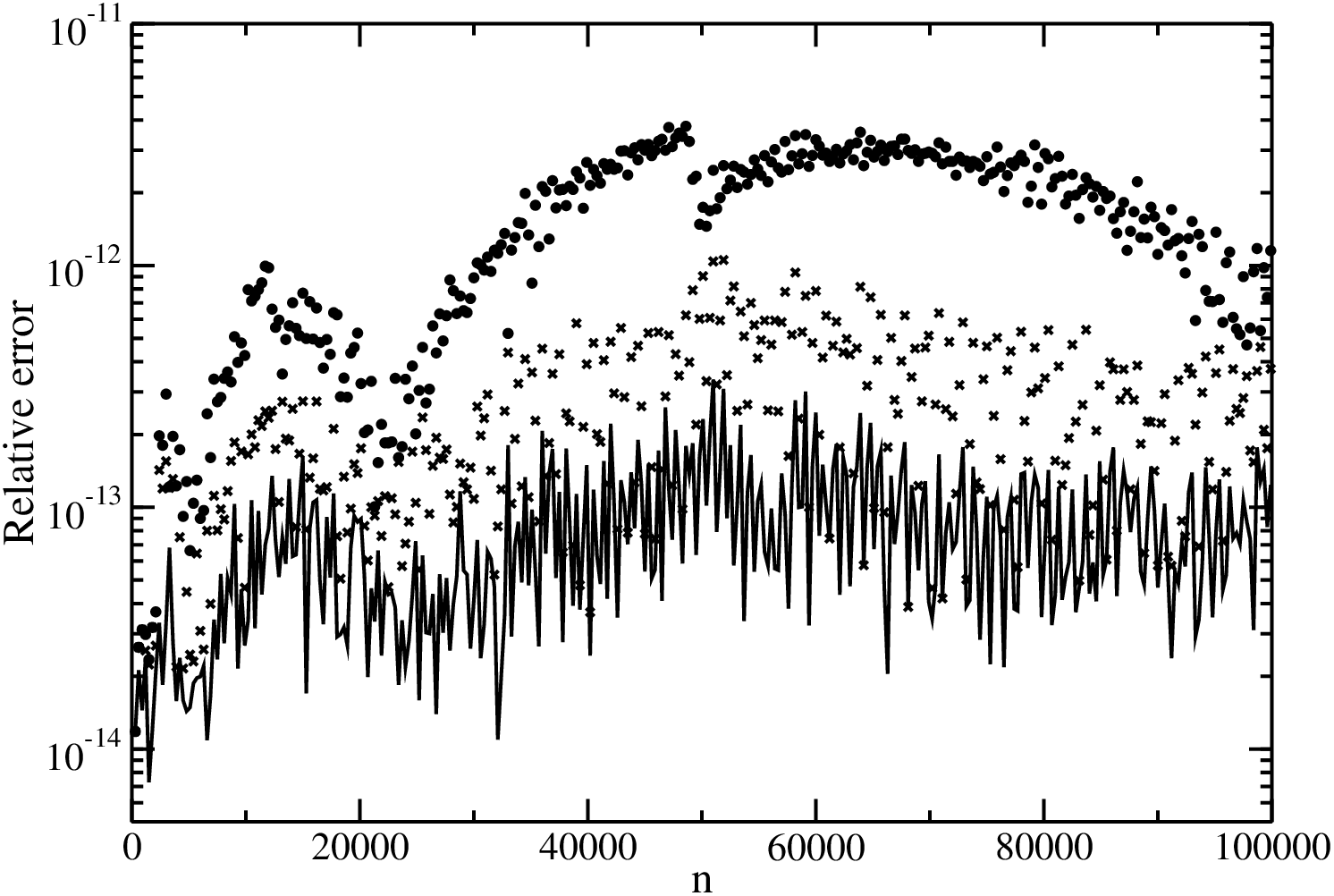}}
\caption{Relative errors in the computation of the weights for $n$-point Gauss--Hermite
quadrature. The dots represent the values 
$\max|1-\omega_{i}^{(d)}/\omega_i^{(q)}|$, where $\omega_i^{(d)}$ are scaled 
weights computed in double precision and 
$\omega_i^{(q)}$ are the weights in quadruple precision. The crosses and the solid line represent
the maximum error when it is evaluated only for the nodes for which the (unscaled) weights
$w_i$ are larger than $10^{-300}$ (crosses) or $10^{-30}$ (solid line).
}
\label{fig2}
\end{figure}

For the scaled weights, the error degradation is moderate, as Figure~\ref{fig2} shows.
The largest errors always correspond to the weights $\omega_i$ for the largest zeros, that is,
to the least significant unscaled weights. If we only compute the errors corresponding
to nodes for which the unscaled weights $w_i$ are larger than a given threshold (say $10^{-30}$ or
$10^{-300}$)
then the errors can be reduced, as shown in Figure~\ref{fig2}.

\begin{figure}
\vspace*{0.8cm}
\centerline{\includegraphics[height=6cm,width=12cm]{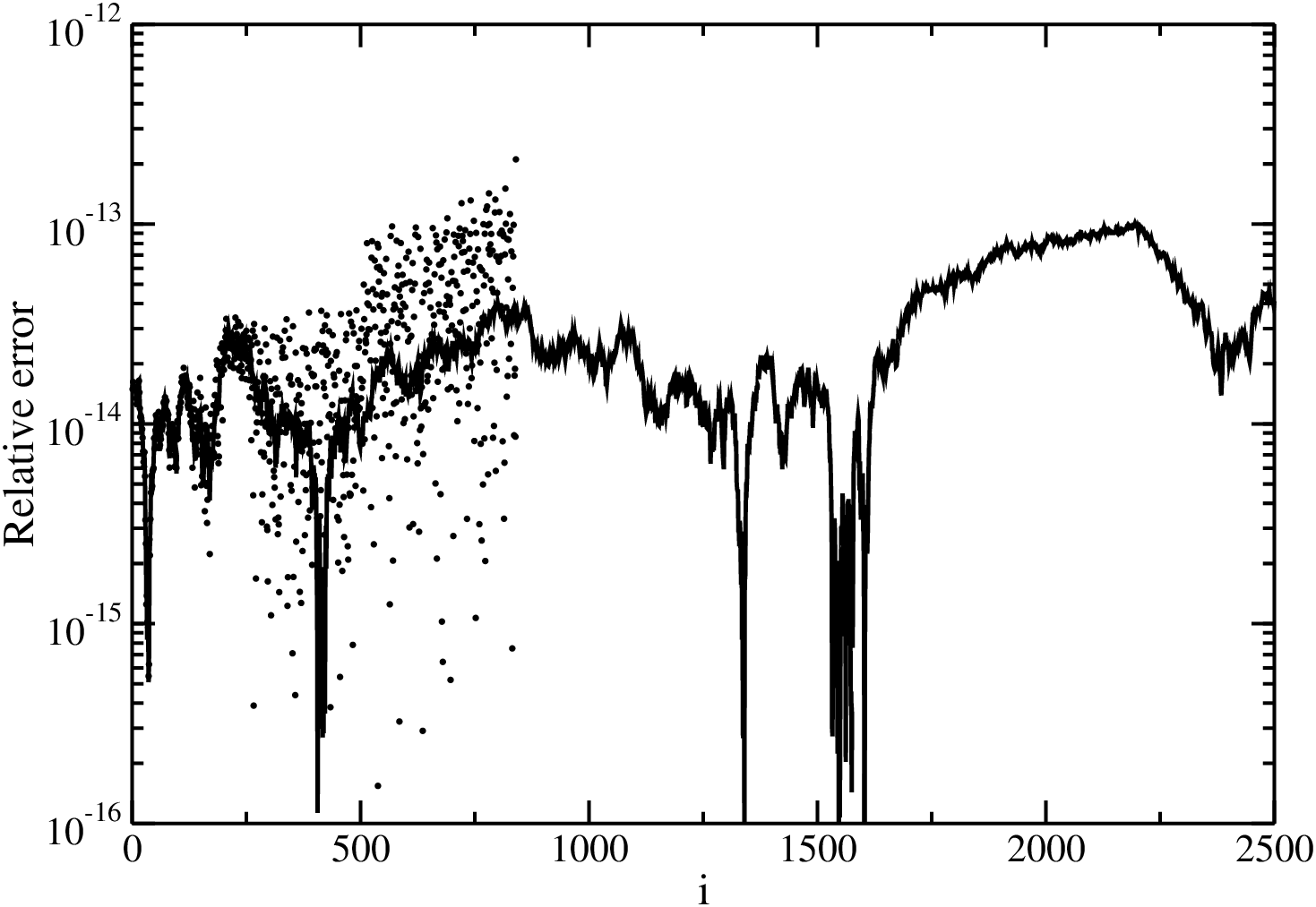}}
\caption{Relative errors in the computation of the weights for $5000$-point Gauss--Hermite
quadrature as a function of $i$, with $i$ numbering the positive nodes in 
increasing order. The solid line corresponds to the errors of the scaled weights and the dots correspond to the unscaled weights. Only the errors for unscaled weights which are larger than the underflow limit are shown.
}
\label{fig3}
\end{figure}

For computing the unscaled weights $w_i$, we have to multiply by the exponential
factor, which gives an additional error, as shown in Figure~\ref{fig3}. Of course, this
figure  shows the relative errors only for those unscaled weights which are larger than the double precision
underflow limit (roughly $25\%$ of all the weights).

As we have shown, the method is more accurate than previous methods, particularly
for the nodes but also for the weights. The method is also fast, and has clear
advantages with respect to the GLR algorithm in terms of complexity. In the
first place, the use of a Runge--Kutta method for computing a first approximation
to the nodes is not needed, because the method generates in all instances monotonic
convergent sequences and automatically provides first estimations which are more
accurate as the order increases. In the second place, the fixed point method is
of order $4$, and the number of iterations for each zero can be made smaller 
by using this fact; then, for instance, if the goal is to compute a node
with a relative accuracy $10^{-16}$, we can stop the iteration safely when two
consecutive iterations satisfy $|x^{(k+1)}-x^{(k)}|<6^{1/4} 10^{-4}$, because the
fact that the order is $4$ implies that the relative error for $x_{k+1}$ can be estimated to be close
to $10^{-16}$, which is smaller than the double 
precision machine-epsilon\footnote{This is a consequence of the absolute error relation (see \cite[Eq.~(2.13)]{Segura:2010:RCO})
$x^{(k+1)}-\alpha \approx \Frac{1}{12}A'(\alpha)(x^{(k)}-\alpha)^4$, with $\alpha$ the root that
is computed, together with the fact that for the Hermite equation $A(x)=2n+1-x^2$, which implies that 
$1-\Frac{x^{(k+1)}}{\alpha} \approx \Frac{1}{6}(x^{(k)}-\alpha)^4\approx \Frac{1}{6}(x^{(k)}-x^{(k+1)})^4$}. With this, we only require a few iterations, and fewer iterations
are required as the order is larger. For instance, for $n=100$ we require $1$ or $2$
iterations per root, and only one for $n>1000$. For orders smaller than $100$ only a few
nodes require $3$ iterations.

\begin{figure}
\vspace*{0.8cm}
\centerline{\includegraphics[height=6cm,width=12cm]{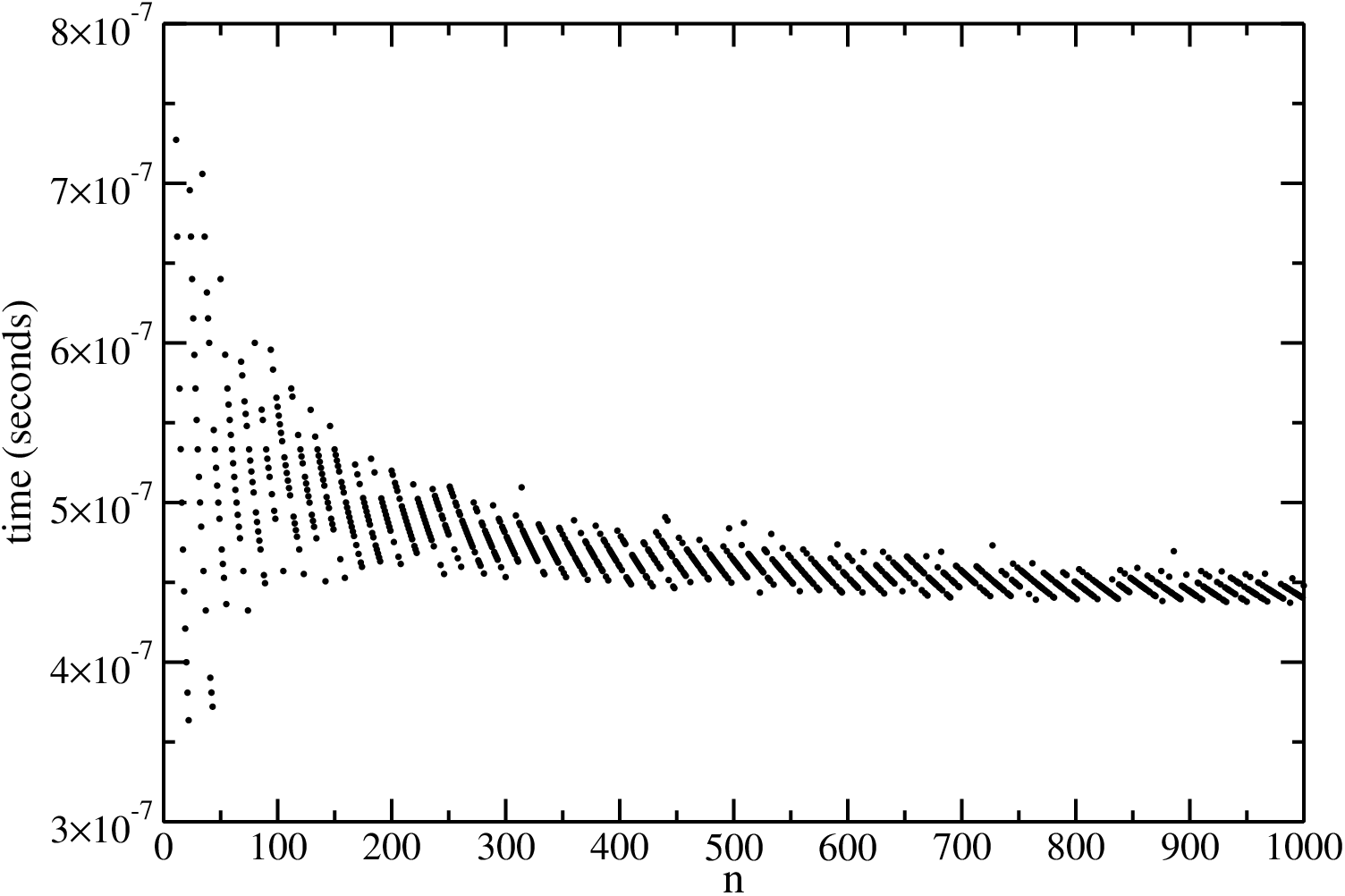}}
\caption{Unitary CPU-time spent (in seconds) as a function of the degree $n$ for Gauss--Hermite quadrature.}
\label{figCH}
\end{figure}

Figure~\ref{figCH} 
shows the CPU times for the computation of Gauss--Hermite quadratures of degree
smaller than $1000$. We show the CPU time divided by $n$, and therefore the figure shows
the CPU-time spent on each node and its corresponding weight. We observe that
this unitary time decreases moderately as $n$ increases approaching an
asymptote, as expected. 
 
The natural
comparison of our method is with the method of \cite{Glaser:2007:AFA}, which is also a purely
iterative method with no asymptotics involved (at least for the Hermite case), and with respect to that
method we have the important advantage that convergence is certain and that our iterative
method doubles the order of convergence of the Newton method considered in that paper.
We conclude that our method should be faster than the one given in \cite{Glaser:2007:AFA} 
and, in any case, it is more accurate. From the comparison of our method with  
\cite{Townsend:2016:FCO} (based on asymptotics and the Newton method) we also conclude that
our method is preferable in terms of accuracy, because that method only provided absolute 
accuracy for the small nodes, and we obtain full accuracy for all the nodes. 
The method in \cite{Townsend:2016:FCO} is based on asymptotics, both for providing first 
estimations to the nodes and for computing the polynomials when applying the Newton iteration.

Recently, and similarly to what was done in \cite{Bogaert:2014:IFC} for Gauss---Legendre, we
provided purely asymptotic methods for computing Gauss--Hermite (and also Gauss--Laguerre)
quadratures. As discussed in that paper, the accurate computation of the nodes 
(and weights) with asymptotics is only slightly more expensive than computing the simpler 
estimations in \cite{Townsend:2016:FCO}, but we have the additional advantage that the Newton
iteration is skipped, therefore speeding up the method. 
This should provide one the fastest method 
of computation for moderately large degree, similarly as \cite{Bogaert:2014:IFC} is the 
fastest method for Gauss--Legendre quadrature of moderately 
large degree.\footnote{A fair comparison of efficiency between different methods should be always 
made by using implementations in a same platform 
and programming language (as is the case of the codes in this paper and \cite{Gil:2018:GHL}). 
This means that the different codes (if available) 
should be translated to a same language and carefully optimized. 
It would be certainly interesting
to coherently benchmark the different available approaches, but 
this is outside the scope of the present paper.}
But the purely iterative method presented in this paper is so efficient that it is even slightly 
faster than the direct computation by 
asymptotics (without iterations) given in \cite{Gil:2018:GHL} (compare Figure~\ref{figCH} with Table 1 in
that reference). 

However, the asymptotic methods in \cite{Gil:2018:GHL} are more accurate for the
computation of the weights in fixed double precision. But the present method has the advantage that 
it works for arbitrary
precision, being based on convergent processes, and that it is valid for any degree and not only for large
degree. In addition, it has the advantage over all the rest of methods that it is a method of fast convergence
and without practical restrictions in the degree.
It is also worth noting the extreme simplicity of the resulting method: in our implementation, 
only around 100 code lines are needed. 
 
We provide several draft codes implementing the methods described in 
this paper\footnote{See http://personales.unican.es/segurajj/gaussian.html}. In particular, for 
the Gauss--Hermite case we provide a Maple worksheet and two Fortran 95 codes: one for double and another
one for quadruple precision. The Maple worksheet can be used for very high accuracy computations, and we
have tested the algorithms for computations with more that $1000$ correct digits. At the end of this section, we
 discuss in some 
detail these high accuracy computations.
 
\subsubsection{A finite precision subtlety when computing Taylor series}
\label{subtle}
We notice that the algorithm is able to produce the Hermite nodes with full double precision and without error
degradation. Several features may explain this fact. The fact that we are computing zeros in increasing order is
favourable for the stability, and a second important fact is that the function values for starting the process are exact 
and without rounding errors (see the discussion just before Eq.~(\ref{LTS})); 
this, together with the fact that in this case the Taylor series have infinite radius of convergence, are factors
which contribute to the stability of the method. Much care must be taken in the implementation of the algorithms
in finite precision arithmetic in order to exploit all these good properties. We mention here a subtle 
programming detail which may result in accuracy loss if not correctly taken into account.

In our method, given an iterate $x^{(i)}$ we can compute the next iterate by taking 
$x^{(i+1)}=T_{-1}(x^{(i)})$ (assuming that $A(x)$ is decreasing, as is the case for Gauss--Hermite when $x>0$);
let us write $T_{-1}(x^{(i)})=x^{(i)}+\delta (x^{(i)})$, where 
$$\delta (x)=-\Frac{1}{\sqrt{A(x)}}\arctan_{-1} (
\sqrt{A(x)}y(x)/y'(x)).$$ As explained before, after the new iteration $x^{(i+1)}$ has been computed, the values
of $y(x^{(i+1)})$ and $y'(x^{(i+1)})$ needed to compute the next iteration are evaluated by using Taylor series centered at $x^{(i)}$ with step
$\delta (x^{(i)})$ (and the process is repeated until an accurate approximation to the node is obtained). 
Two different ways for computing $y(x^{(i+1)})$ and $y'(x^{(i+1)})$ with Taylor series are the following:

\begin{enumerate}
\item{}$h=\delta (x^{(i)})$; $x=x^{(i)}$;  $x^{(i+1)}=x^{(i)}+h$;
\item{}Use the Taylor series (\ref{LTS}), with derivatives computed by (\ref{reH}).
\end{enumerate}

\begin{enumerate}
\item{}$x^{(i+1)}=x^{(i)}+\delta (x^{(i)})$; $x=x^{(i)}$; $h=x^{(i+1)}-x^{(i)}$;
\item{}Use the Taylor series (\ref{LTS}), with derivatives computed by (\ref{reH}).
\end{enumerate}

Notice that we have just interchanged the order of the first and the third evaluations, 
but this results in two noticeable different ways to compute the series in 
finite precision arithmetic. 
Initially, one could think that the 
second option is worse, because we are computing a small quantity ($h$) as the difference of two quantities
($x^{(i)}$ and $x^{(i+1)}$) which are typically  much larger, and it is common wisdom that this will
 introduce rounding errors in $h$ (and then in the computation of Taylor series). 
But the situation here is the opposite, and the second option is preferable. 

Let us concentrate in the first implementation. In that first option we start by computing $h$, which we
can do accurately; then we take $x=x^{(i)}$ and we compute  $x^{(i+1)}=x^{(i)}+h$. 
We notice that, because $h$ can be much smaller than 
$x$, the information carried by a number of the digits of $h$ will be lost when computing $x^{(i+1)}=x^{(i)}+h$. 
Because of this, in finite precision we have that $x^{(i+1)}-x^{(i)}\neq h$, and therefore
 using $h$ for computing the Taylor series is not a good idea: $h$ does not measure faithfully the 
difference $x^{(i+1)}-x^{(i)}$. This is why the first option is the wrong choice. The right choice is the
second one because $h=x^{(i+1)}-x^{(i)}$ does faithfully represent the difference between
successive iterates and therefore it is a better step for the Taylor series.

In our Fortran programs, we use the second option and in this way we avoid all error degradation 
for computing the 
Gauss--Hermite 
nodes; contrarily, when we choose the first option, the algorithms tend to accumulate errors as we compute successive nodes.

\subsubsection{Very high accuracy computations}

The reliability
and high order of convergence makes our iterative algorithms specially suited for high accuracy. 
In order to test the performance of our algorithms for high accuracy computations we have translated our
Fortran Gauss--Hermite quadrature program to Maple, which allows us to test
the methods for very high accuracies (we have tested the algorithm down to $10^{-1024}$ relative accuracy). 
The asymptotic exactness
of the methods implies that the computation per node improves as the degree increases because the method
tends to be exact, and then the number of required iterations decreases; but also, as
we will see, the computation by Taylor series becomes more efficient as the degree increases.  

Table \ref{table} shows the average number of iterations per node together with the number of terms
of the Taylor series per node (summing the total number of the terms in all the iterations needed to
compute the node), both as a function of the degree and the relative accuracy. 
For each node, the iterations are stopped when two consecutive estimations $x_n$ and $x_{n+1}$ 
are such that $|x_{n+1}-x_n|<6^{1/4} 10^{-D/4}$ where
$D$ are the digits of accuracy (in the table, $D=2^{3+E}$); this implies, 
as discussed in Section~\ref{numres}, 
that the error relative error of $x_n$ will be
approximately $10^{-D}$. As for the stopping rule for the series, each sum is terminated when the last computed term gives a relative contribution smaller than $10^{-D}$; we supplement this by requiring that at least $20$ terms are considered in each iteration, and we also limit the maximum number of terms per iteration to
$50\times E^{1.5}$ (this last condition is convenient for moderate degrees).

The first feature
to notice from Table \ref{table} is that, for accuracies up to $1024$ digits, and for degrees up to $10^{5}$, the average number of iterations is never greater than $5$ and that, as expected, this 
number decreases as the degree increases. We observe also that the number of terms for the series
is smaller as the degree increases. On the other hand, as expected, the number of iterations and of terms in the
series increases with the demanded accuracy. 

Compared to the more standard method of computation of Hermite polynomials by using the three-term 
recurrence relation, which needs $n-1$ iterations for computing $H_{n}(x)$ starting from $H_0(x)=1$ and $H_1(x)=2x$,
we observe that the series are more efficient than the recurrence relation for $n>1000$, even for high accuracy, while the
recurrence relation may be interesting for low degrees and high accuracy. We note, however, that the series are not only
useful from the point of view of accuracy, but also because it permits computations which are free of overflows/underflows (both for the Hermite and Laguerre quadratures) and, in addition, they appear to be more stable 
in particular for the Laguerre case (see Section~\ref{TTRR}).

\begin{figure}[tb]
\vspace*{0.8cm}
\begin{center}
\begin{minipage}{3cm}
\centerline{\includegraphics[height=6cm,width=6cm]{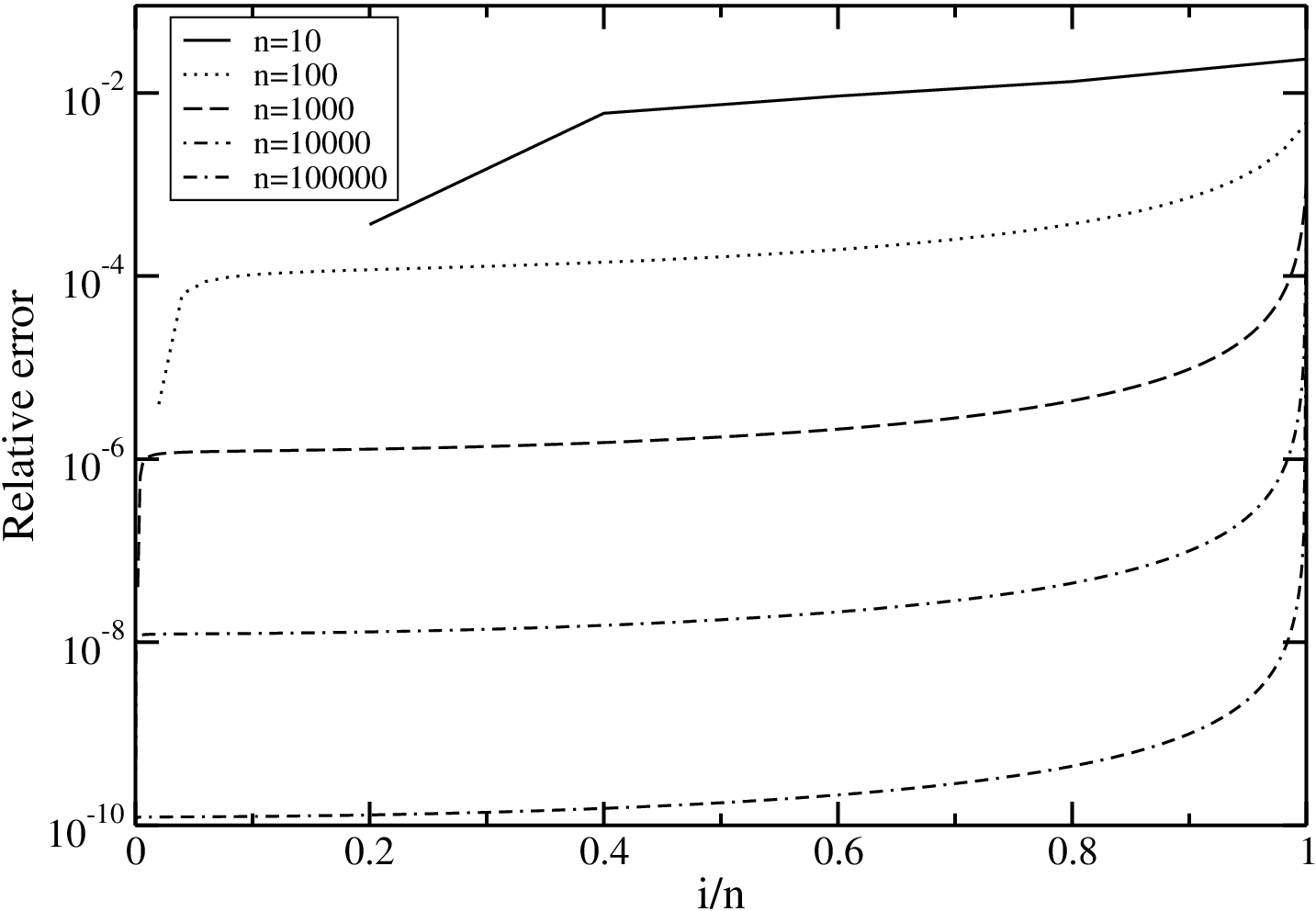}}
\end{minipage}
\hspace*{3cm}
\begin{minipage}{3cm}
\centerline{\includegraphics[height=6cm,width=6cm]{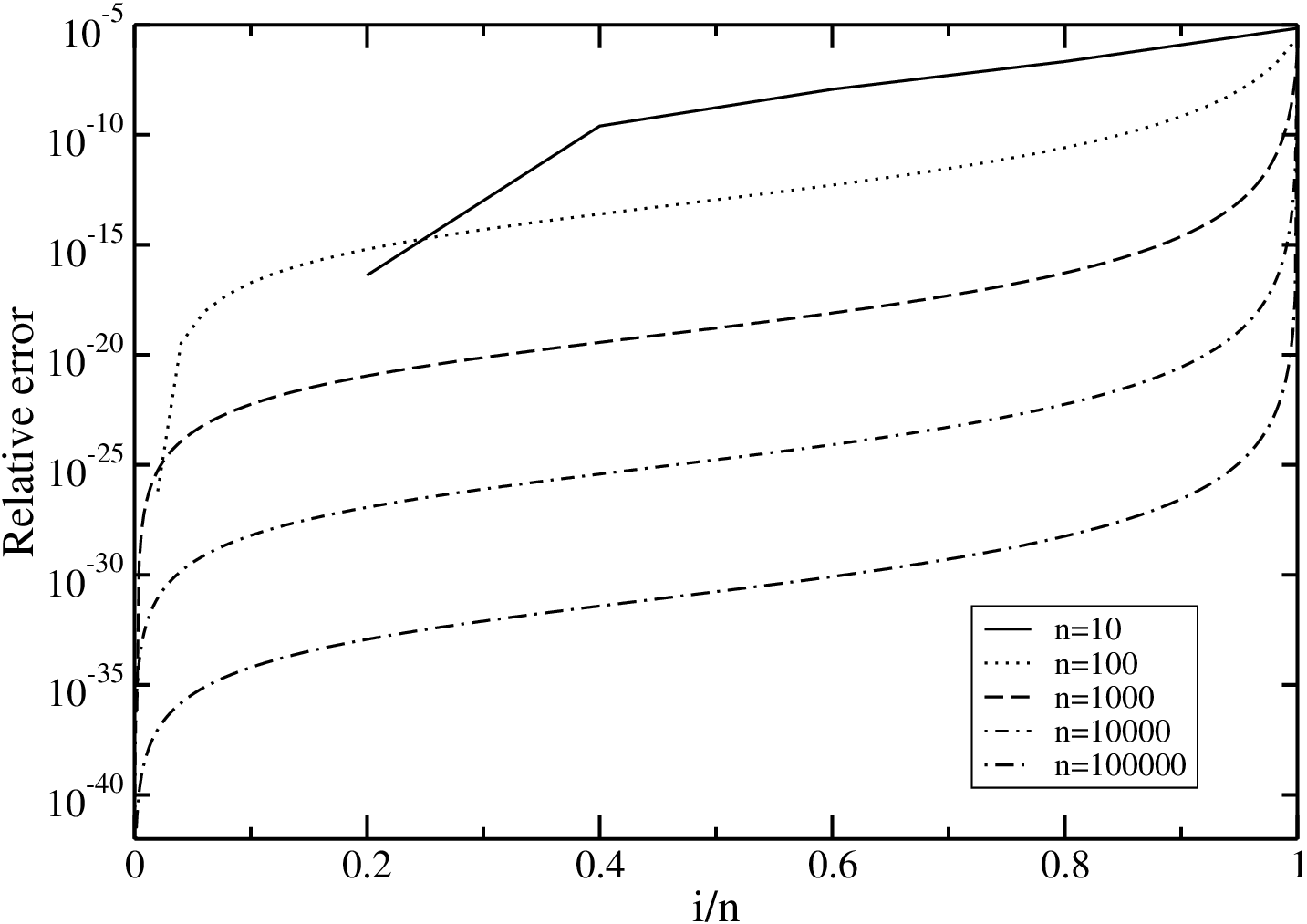}}
\end{minipage}
\end{center}
\caption{Left: relative error in the first estimation for the nodes provided by the algorithm as a functions of 
$i/n$, where $n$ is the degree and $i$ is the index enumerating the positive nodes in increasing order, 
$i=1,\ldots n/2$. Right: same but for the second estimation.}
\label{fig10}
\end{figure}

As explained before, the method is asymptotically exact as the degree tends to infinity in the sense that it
is becomes exact in this limit and, given a node, the next node can be computed in just one step if the degree
is high enough. This fact is illustrated in Figure~\ref{fig10}, and in particular in the figure on the left,
which gives the accuracy for the first estimation of each node; we observe that, as the degree 
increases, the first
estimation becomes more accurate. The figure on the right, on the other hand, 
shows the relative error of the second estimation
and shows the rapid convergence of the method, very specially for the first zeros.

\begin{table}
$$
\begin{tabular}{c|ccccccc}
n \textbackslash E & 1 & 2 & 3 & 4 & 5 & 6 & 7 \\
\hline
 10    & 2/75 & 2.6/117 & 2.8/197 & 3.6/384 & 3.8/747 & 4.6/1522 & 4.8/3128\\ 
 100   & 1.9/69 & 2.1/91 & 2.9/142 & 3.1/222 & 3.9/377 & 4.1/673 & 4.9/1226\\ 
 1000   & 1.2/54 & 2/87 & 2.3/121 & 3/192 & 3.3/308 &  4/528 & 4.3/923 \\
10000   & 1/49 & 2/88 & 2/113 & 3/183 & 3/286 & 4/490 & 4/841 \\
100000 &  1/48& 1.8/82 & 2/112 & 2.9/177 & 3/276 & 3.9/469 & 4/802\\
\end{tabular}
$$
\caption{Average number of iterations per node and number of the terms of the Taylor series used 
per node for degrees 
$n=10,\,100,\,1000,\,10000,\,100000$ (rows) and for relative accuracies of $10^{-D}$ with $D=2^{3+E}$
 (columns).}%
\label{table}%
\end{table}

\section{Computation of Gauss--Laguerre quadratures}
\label{Gauss--Laguerre}

The methods for Gauss--Laguerre quadrature that we next describe,
 as happened for the Hermite case, work for any prescribed
accuracy and for practically unrestricted values of the parameters. 

Here we will describe algorithms corresponding to the Liouville transformation with
the change of variables $z=\sqrt{x}$. As explained before, this is the natural selection
in the sense that the algorithm is asymptotically exact as $n\rightarrow +\infty$ and also
because the nodes are computed in decreasing order of significance of the weights.

The main method of computation of the function $y(z)$ of (\ref{yzl}), or a scaled version, 
is the use of local Taylor series. As in the Hermite case, we use Taylor series for a
conveniently normalized function so that overflows/underflows are avoided, and later
rescale the weights using one of the moments. However, differently from Hermite and for the
reasons above explained, it is not possible to use Taylor series in all occasions, and for
a few nodes/weights we will need alternative methods.

In order to avoid overflows/underflows, these alternative methods (recurrence relation
and a continued fraction) will use ratios
of Laguerre polynomials instead of the polynomials. These alternative methods are needed 
for the first node and in some cases for an additional node when $\alpha$ is small.

It is important to note that the function normalization
 must be consistently maintained for the computation of all the weights, which is guaranteed if the
Taylor series method is used serially in all steps and without normalization changes, 
with the exception of the first computed scaled weight 
which can be chosen arbitrarily. 
This means that 
we can compute the first node and 
the scaled weight with an alternative method and then compute the rest with Taylor series. 
But when an additional node/weight
needs to be computed  with methods different from Taylor series, they must  
be recomputed later by using Taylor series in order to ensure a 
consistent normalization, as we later explain.

\subsection{Computing the functions}

Next we describe the methods of function computation, starting with the two ``exceptional"
methods (recurrence relation and continued fraction) and continuing with the core method (Taylor series).

\subsubsection{Recurrence relation over the degree}
\label{TTRR}
As any orthogonal polynomial, Laguerre polynomials satisfy a three-term recurrence relation
that we can use for computing polynomial values. We write this in terms of ratios.
$$
\begin{array}{@{}r@{\;}c@{\;}l@{}}
R_k (x)&=&L_{k+1}^{(\alpha)}(x)/L_k^{(\alpha)}(x),\\[8pt]
R_0 (x)&=&\alpha+1-x,\\
R_{k} (x)&=&\Frac{1}{k+1}\left[(2k+\alpha+1-x)-\Frac{k+\alpha}{R_{k-1}(x)}\right].\\[8pt]
\end{array}
$$
Then, using the differential relation 
$x {L_n^{(\alpha)}}^{\prime}(x)=n L_n^{(\alpha )}(x)-(n+\alpha)L_{n-1}^{(\alpha)}(x)$ and
the definition
Eq.~(\ref{yzl}) we get
\begin{equation}
\Frac{\dot{y}(z)}{y(z)}=\Frac{2n+\alpha+1/2}{z}-z-\Frac{2(n+\alpha)}{zR_{n-1}(z^2)}.
\end{equation}

This three-term recurrence relation for Laguerre polynomials is not badly conditioned. However,
it should be considered for not too high degrees for two reasons: firstly, because the computation
is not efficient, and secondly because there is some degradation when large orders are considered.
We use the recursion only when $n<10$.

An interesting alternative with fast convergence (if $x$ is not too large) and which is more reliable than the previous recurrence
relation is given next
in terms of a continued fraction.

\subsubsection{Continued fraction}

Considering the relation between the Laguerre polynomials and the Kummer function
$$
L_n^{(\alpha)}(x)=\Frac{(\alpha+1)_n}{n!}M(-n,\alpha+1,x),
$$
and because the Kummer function $M$ satisfies a three-term recurrence relation relating
three consecutive values of $\alpha$ and $M$ is a minimal solution of this recurrence relation
as $\alpha\rightarrow +\infty$ \cite{Segura:2008:NSS}, we deduce that $L_n^{(\alpha)}(x)$ 
is minimal with respect to recursion over $\alpha$ as $\alpha\rightarrow +\infty$.
Therefore $L_n^{(\alpha)}(x)$ is a minimal solution  of the 
recurrence relation
$$
\begin{array}{l}
L_n^{(\alpha +1)}(x)+b_{\alpha}L_n^{(\alpha)}(x)-a_{\alpha}L_n^{(\alpha-1)}(x)=0,\\
b_{\alpha}=-(1+\alpha/x),\quad a_{\alpha}=-(n+\alpha)/x,
\end{array}
$$
as $\alpha \rightarrow +\infty$.
In terms of the ratios $r^{(\alpha)}=L_n^{(\alpha)}(x)/L_n^{(\alpha -1)}(x)$ we have
\begin{equation}
\label{reco}
r^{(\alpha)}=\Frac{a_{\alpha}}{b_{\alpha}+r^{(\alpha+1)}},
\end{equation}
and Pincherle's theorem guarantees that the continued fraction resulting from the iteration
of (\ref{reco}) is convergent, and that it converges to 
$L_n^{(\alpha)}(x)/L_n^{(\alpha -1)}(x)$. Therefore,
$$
r^{(\alpha)}(x)=\Frac{L_n^{(\alpha)}(x)}{L_n^{(\alpha -1)}(x)}=\Frac{a_\alpha }{b_{\alpha}+}
\Frac{a_{\alpha+1}}{b_{\alpha+1}+}\cdots,
$$
and using the derivative rule
$$
x{L_n^{(\alpha)}}^{\prime}(x)=-\alpha L_n^{(\alpha )}(x)+(\alpha+n)L_n^{(\alpha-1)}(x)
$$
and the definition
Eq.~(\ref{yzl}) we obtain
\begin{equation}
\label{CFcoc}
\Frac{\dot{y}(z)}{y(z)}=\Frac{1/2-\alpha}{z}-z+\Frac{2(n+\alpha)}{zr^{(\alpha)}(z^2)}.
\end{equation}

\subsubsection{Computation of Taylor series}

The function $y(z)$ satisfies (see (\ref{normL})) 
\begin{equation}
\label{ecudef}
P(z)y^{(2)}(z)+Q(z)y(z)=0,
\end{equation}
with
$$
P(z)=z^2,\quad Q(z)=-z^4+2 L z^2+\Frac{1}{4}-\alpha^2.
$$
Taking successive derivatives and using that $P^{(n)}(z)=0$, $n>2$ and $Q^{(n)}(z)=0$, $n>4$, we
obtain the following recursion formula for the derivatives with respect to $z$:
\begin{equation}
\label{taylag}
\displaystyle\sum_{m=0}^{2}
\left(
\begin{array}{c}
j\\
m
\end{array}
\right)
P^{(m)}(z)y^{(j+2-m)}(z)+\displaystyle\sum_{m=0}^{4}
\left(
\begin{array}{c}
j\\
m
\end{array}
\right)
Q^{(m)}(z)y^{(j-m)}(z)=0,
\end{equation}
where $\left(
\begin{array}{c}
j\\
m
\end{array}
\right)$
are binomial coefficients. 

Eq.~(\ref{taylag}) is a seven-term recurrence relation ($y^{(j)}$ appears in the last term of the first sum and in the first term of the second sum), and therefore the space of solutions has dimension 6.
Considering the Perron-Kreuser theorem \cite{Kreuser:1914:UDV}, the solutions of this difference
equation lie in two subspaces: a subspace of dimension two of solutions satisfying
$$
\limsup_{n\rightarrow +\infty} |y^{(n)}/n!|^{1/n}=|1/x|,
$$
and a subspace of dimension four of solutions satisfying
$$
\limsup_{n\rightarrow +\infty} |y^{(n)}/\sqrt{n!}|^{1/n}=1.
$$
The solutions of the first subspace are dominant over the second subspace. The 
derivatives of solutions of (\ref{ecudef}) are in this dominant subspace, and
it contains functions
which have Taylor series centered at $x$ of radius $R=|x|$ (as corresponds to a differential equation
with a singularity at $x=0$). Because of the dominance of these solutions, the computation of the
derivatives in the forward direction is well conditioned.

Of course, this is not the only Taylor series that we could consider, and we could use series for other
functions, and also in other variables; for instance, we could consider a Taylor series for 
(\ref{norlx}), which satisfies (\ref{Alx}). However, there are good reasons to use this 
form of Taylor series. Firstly, as we commented before,
as $n\rightarrow +\infty$, the ODE is such that the coefficient is essentially constant in the largest part of the
interval of oscillation. This means that the solutions will have a slowly varying amplitude of 
oscillation (and also a slowly varying period of oscillation), and this reduces drastically
the possibility of overflows/underflows in the computation. In the second place, as discussed
before, the conditioning of the computation of the scaled weights
 is very good, as they do not
depend at first order approximation on the value of the nodes.

\subsection{Computation of the nodes}

As we discussed in Section~\ref{Lagcha}, when we consider the Liouville transformation with 
change $z(x)=\sqrt{x}$ we have to distinguish between two cases: $|\alpha|\le 1/2$
and $|\alpha|> 1/2$. 

In the algorithms it will be useful to consider the bounds given in \cite{Dimitrov:2010:SBF}, which
can be written as:
\begin{equation}
\begin{array}{l}
x_u=\Frac{2 n^2+n (\alpha-1)+2(\alpha+1)+2(n-1)\sqrt{n^2+(n+2)(\alpha+1)}}{n+2},\\
\\
x_l=P/x_u,\quad P=\Frac{(\alpha+1)(n(\alpha+5)+2(\alpha-1))}{n+2}.
\end{array}
\end{equation}
All the zeros of $L_n^{(\alpha )}(x)$, $\alpha>-1$, are in the interval $(x_l,x_u)$ and we observe that,
for large $n$, $x_u=4n(1+{\cal O}(n^{-1}))$ and $x_l=\frac14(\alpha+1)(\alpha+5)n^{-1}(1+{\cal O}(n^{-1}))$.

For the first case, $|\alpha|\le 1/2$, because the coefficient $A(z(x))$ of
(\ref{normL}) is decreasing for positive $x$,
we can start the algorithm at the lower bounds
$z=z_l=\sqrt{x_l}$, and compute
the zeros in increasing order until all the nodes are computed.

For the second case we would start at $z_e=(\alpha^2-1/4)^{1/4}$ and 
then compute the zeros larger than $z_e$ in increasing order with $T_{-1}$ (see (\ref{def1}))
and the smaller zeros in decreasing order with $T_{+1}$. 
When computing the larger zeros in increasing order we can
stop the computation when the upper bound for the zeros is
surpassed.\footnote{Note that all the $z$-values generated by $T_{-1}$ 
are an increasing sequence and therefore the stopping rule
is safe.}  Then, the smaller zeros (if any) are computed until the total number of nodes is completed.

Let us observe that the coefficient $A(z(x))$ (\ref{normL}) is positive in the
interval $(x_L,x_R)$, with $x_R=L+\sqrt{L^2-\alpha^2+1/4}$ and 
$x_L=(\alpha^2-1/4)/x_R$. One can prove that $x_R>x_u$ for all $\alpha>-1$
and $x_L<x_l$ if $\alpha>-7/8$. This means that $(x_l,x_u)\subset (x_L,x_R)$ if
$\alpha>-7/8$, and therefore $A(z(x))$ is positive in an interval 
containing all the nodes. However, when $\alpha<-7/8$ the smallest zero $x_1$
can be such that $A(z(x_1))<0$, and this is certainly so as $\alpha\rightarrow -1$
because the first node $x_1$ tends to zero in this limit. Observe that this may
only happen for the smallest zero, since only one zero may exist in the
interval $(0,x_L]$ (this is simple to check by analyzing the monotonicity/convexity 
of the solutions of the ODE, as done, for instance, in 
\cite[Sect. 3.1]{Segura:2010:RCO}). 

In \cite[Sect. 3.1]{Segura:2010:RCO} it is discussed how to modify the 
fixed point method for a reliable fourth-order convergence to the first node $x_1$
when $A(z(x_1))<0$. In our case, we are working with the variable $z=\sqrt{x}$,
and the modification consists in applying the fixed
 point method $z^{(n+1)}=T(z^{(n)})$
\begin{equation}
\label{modi}
T(z)=z-\Frac{1}{\sqrt{-A(z)}}\mbox{atanh}\left(\sqrt{-A(z)}h(z)\right)
,\quad h(z)=y(z)/\dot{y}(z),
\end{equation}
if $A(z^{(n)})<0$, instead of the fixed point methods 
$$
T_{\pm 1}(z)=z-\Frac{1}{\sqrt{A(z)}}\arctan_{\pm 1}\left(\sqrt{A(z)}h(z)\right),
$$
which we use when $A(z)>0$. 
We notice that when $|\alpha|>1/2$ we always start from a value 
(also when $A(z(x_1))<0$) $x_e = \sqrt{\alpha^2-1/4}$ such that $A(z(x_e))>0$, 
which means that,
as described above, we can safely start the algorithm by using $T_{-1}$
to compute the zeros larger than $x_e$ and then use 
$T_{+1}$ to compute those smaller than $x_e$,
switching to (\ref{modi}) for the smallest zero if needed. This scheme converges
with certainty. 

With respect to the methods of computation of Laguerre polynomials (or related functions), 
we can also distinguish between two cases.
The simplest case is when $|\alpha|>1/2$ and $n$ is sufficiently large. 
Observe that for fixed $|\alpha|>1/2$
the number of nodes smaller than $z_{e}=z(x_e)$ is ${\cal O}(\sqrt{n})$ as 
$n$ increases.\footnote{We have $z_e-z_{l}=(\alpha^2-1/4)^{1/4}(1+{\cal O}(n^{-1/2}))$ and the
maximum value of $A(z)$ is reached at $z=z_e$, where 
$A(z_e)=2(L-\sqrt{\alpha^2-1/4})$. Therefore, the distance between consecutive nodes
in the $z$ variable can be bounded by $\Delta z=z_{i+1}-z_{i}>\pi/\sqrt{A(z_e)}$;
with this we estimated that the number of zeros smaller than $z_e$ is 
$(z_e-z_l)\sqrt{A(z_e)}/\pi\sim \Frac{2(\alpha^2-1/4)^{1/4}}{\pi}\sqrt{n}$, which is an upper bound.}
This implies that for large enough $n$ the computation of all the nodes can
be carried out using only Taylor series, except for computing the starting value at $z=z_e$.

Let us recall that the use of Taylor series is limited by the fact that $z=0$ is a singularity
of the differential equation, and that the Taylor series centered at a $z>0$ has a radius of convergence 
$R=z$. In the case when there are zeros smaller that $z_e$ ($|\alpha|>1/2$, $n$ large enough), 
we never need to evaluate Taylor series outside its radius of convergence, because
we compute in the direction of decreasing $z$; then the use of Taylor
series is safe (also because we do not need to evaluate series very far away from their center).

The situation is different for $|\alpha|\le 1/2$, but also for slightly larger $|\alpha|$ and
$n$ small. Here, not only we need to start with the CF (or recurrence relation) for the first zero, but also
we need the CF for the second zero. In the case $|\alpha|\le 1/2$ we have, because the $A(z)$ 
coefficient (\ref{normL}) is decreasing
and $z=0$ is a zero of $y(z)$ (\ref{yzl}), that $z_1-0<z_2-z_1$, where $z_1<z_2$ are the two 
smallest positive zeros of $y(z)$. This means that the disc of absolute convergence of the series
centered at $z_1$, which has radius $R=z_1$ does not include $z_2$. This indicates that 
Taylor series should not be used for computing $z_2$ after $z_1$ has been computed.

We conclude that for $|\alpha|\le 1/2$ we need the computation of the CF for evaluating the first two
zeros; but also for larger $|\alpha|$ the use of Taylor series may be inaccurate for the first zeros,
particularly for small $n$. For instance, for $n=4$, $\alpha=1$ we
have $(\alpha^2-1/4)^{1/4}=0.930604$, 
$z_1=0.86214380$, $z_2=1.60363182$, and again $z_1>z_2-z_1$.

In the case that both the forward and the backward sweep ($T_{-1}$ and $T_{+1}$) are used ($|\alpha|>1/2$,
 large enough $n$),
the algorithm is in its simplest form, and only one evaluation of the recurrence relation or the CF is required for starting the process, after which Taylor series expansions are used. Also in the case when $n$ is smaller 
and no backward sweep is needed,
we only need one CF evaluation provided that $z_1$ is sufficiently larger than 
$z_2-z_1$, because in this case the Taylor series can be used to go from $z_1$ to $z_2$. 

In practice, when $\alpha\ge 2$ we only require one CF (or recurrence relation)
evaluation, while for $-1<\alpha<2$, although it is not necessarily in all cases,
we prefer to use the CF for the evaluation of the first two nodes larger
than $z_e$. As described above, Taylor series can be safely used for computing
the nodes smaller than $z_e$ by a backward sweep. 

\subsection{Computation of the weights}

The scheme for computing the weights depends on the number of nodes for which
the CF is required. We first describe the simplest case when only one CF is
required.

\subsubsection{With only one CF evaluation}

When all the zeros larger than $z_e$ can be accessed 
with Taylor series (case $\alpha>2$ in our algorithm), the computation
of the weights and scaled weights goes as follows.

We start at $z=z_e$ by computing $h=y/\dot{y}$ with the CF (or the recurrence relation).
Then we set, for instance $\bar{y}(z_e)=h(z_e)$ and 
$\dot{\bar{y}}(z_e)=1$ \footnote{Or, if $h(z_e)$ is very large, we can instead
take $\bar{y}(z_e)=1$ and $\dot{\bar{y}}(z_e)=1/h(z_e)$ to prevent overflows.}, 
where we denote by $\bar{y}$ a solution of the ODE (\ref{normL}). 
Now, we make the rest of the computations using Taylor series,
with initial values given by $\bar{y}(z_e)$ and $\dot{\bar{y}}(z_e)$. 

In the same computation of the nodes, we will obtain numerical approximations
for $\dot{\bar{y}}(z_i)$, and then we obtain the scaled weights $\omega_i$ 
(\ref{wasym2}) up to a factor (say $\gamma$):
$$
\bar{\omega_i}=|\dot{\bar{y}}(z_i)|^{-2}=\gamma |\dot{y}(z_i)|^{-2}=\omega_i .
$$
This factor can be fixed by normalizing with the first momentum, that is,
\begin{equation}
\label{muce}
\mu_0=\displaystyle\sum_{j=0}^n w_i=\int_0^{+\infty}x^{\alpha}e^{-x}dx=
\Gamma (\alpha+1),
\end{equation}
where 
\begin{equation}
\label{relsu}
w_i=\omega_i x_i^{\alpha +1/2}e^{-x_i},\quad x_i=\sqrt{z_i}. 
\end{equation}

We observe that the normalization (\ref{muce}) may result in overflow problems
in floating point arithmetic
when $\alpha$ is large. For this reason, we prefer to compute weights 
$\hat{w}_i$ normalized to one, that is, such that 
$\sum_{i=1}^n \hat{w}_i=1$, in other words, $\hat{w}_i=w_i/\Gamma (\alpha +1)$
with scaled weights $\hat{\omega}_i=\omega_i/\Gamma (\alpha +1)$.

For this purpose we start with the scaled weights $\bar{\omega}_i$ computed by the
algorithm. From these, we compute unnormalized and
unscaled weights considering the factor in (\ref{relsu}), but we do so 
relative to the first node larger than $z_e$ (which either corresponds to the
largest weight or is close to it). Suppose that this weight is the $j$-th, then
we take 
\begin{equation}
\label{sca}
\bar{w}_i=\bar{\omega}_i \exp (F_i),\quad 
F_i=x_j-x_i+(\alpha+1/2)\log\left(\Frac{x_i}{x_j}\right),\quad i=1,\ldots n.
\end{equation}
We can do this in parallel with the computation of the nodes, and we
can decide, using the fact that the weight $\bar{w}_i$ are computed in decreasing
order of magnitude, how many weights/nodes we need. Of course, computing unscaled weights
$\bar{w}_i$ smaller than the underflow number is unnecessary, but we may decide
to compute the nodes and scaled weights. 

In any case, the unscaled unnormalized weights are related to the weights 
$\hat{w}_i$ (normalized to one) by a factor $\lambda$, 
$\bar{w}_i=\lambda \hat{w}_i$, where obviously 
$\lambda=\displaystyle\sum_{i=1}^n \bar{w}_i$. Then, we compute the unscaled
(and normalized to one) weights $\hat{w}_i=\bar{w}_i /\lambda$ and the 
corresponding scaled weights  $\hat{\omega}_i=\bar{\omega}_i /\lambda$.

\subsubsection{With additional CF evaluations and two backward steps}

As explained before, when $-1<\alpha<2$ we choose to use the CF (or recurrence relation)
for the computation of the first two nodes larger than $x_e$ ($z>z_e$) for $|\alpha|>1/2$ or larger
than the lower bound $x_l$ for $|\alpha|\le 1/2$. Let us denote these
two nodes (in the $z$-variable) by $z_k$ and $z_{k+1}$. After $z_{k+1}$ has been computed, 
we can continue 
with the larger nodes  $z_i$, $i>k+1$, with Taylor series, taking as initial values, $\bar{y}(z_{k+1})=0$ and,
for instance, $\dot{\bar{y}}(z_{k+1})=1$ and proceeding with the iteration $T_{-1}$. 
Then, for all the nodes $x_i$, $i\ge k+1$ the normalization for the scaled weights
$\bar{\omega}_i$ is consistent because we have used Taylor series for all of them. 
For the weight $\bar{\omega}_k$ to be consistent 
with the same normalization, we can use Taylor series centered at $z_{k+1}$ 
to compute $\dot{\bar{y}}(z_k)$; this is all that needs to be done if 
$|\alpha|\le 1/2$. In the case $|\alpha|>1/2$, if the number of computed nodes larger than $z_e$ 
does not equal the degree $n$, then we have zeros smaller than $z_e$; in this case, 
we compute $\bar{y}(z_e)$ and 
$\dot{\bar{y}}(z_e)$ using Taylor series centered at $z_k$, and continue with the computation of
the nodes smaller than $z_e$ with the use of the fixed point method $T_{+1}$ and the application of
Taylor series. This completes the computation of the scaled unnormalized weights $\bar{\omega}_i$,
and we can compute the unscaled weights $\hat{w}_i$ (normalized to $1$) 
and the corresponding scaled weights in the same way as before.

\subsection{Numerical results}

We have implemented our algorithm in a double precision
Fortran routine and we have compared it with a quadruple precision version
of our algorithm. Additionally, we have tested the algorithms against a Maple
implementation of our methods (with Laguerre polynomials computed by Maple commands) 
in order to ensure the correctness of the method.
Differently from the Hermite case, we observe error degradation for the nodes, 
and a moderate error degradation for the weights is also observed. The
source of error comes from the initial value given by the continued fraction
or the recurrence relation, and the propagation of the errors in the application of
Taylor series.

\begin{figure}
\vspace*{0.8cm}
\centerline{\includegraphics[height=6cm,width=12cm]{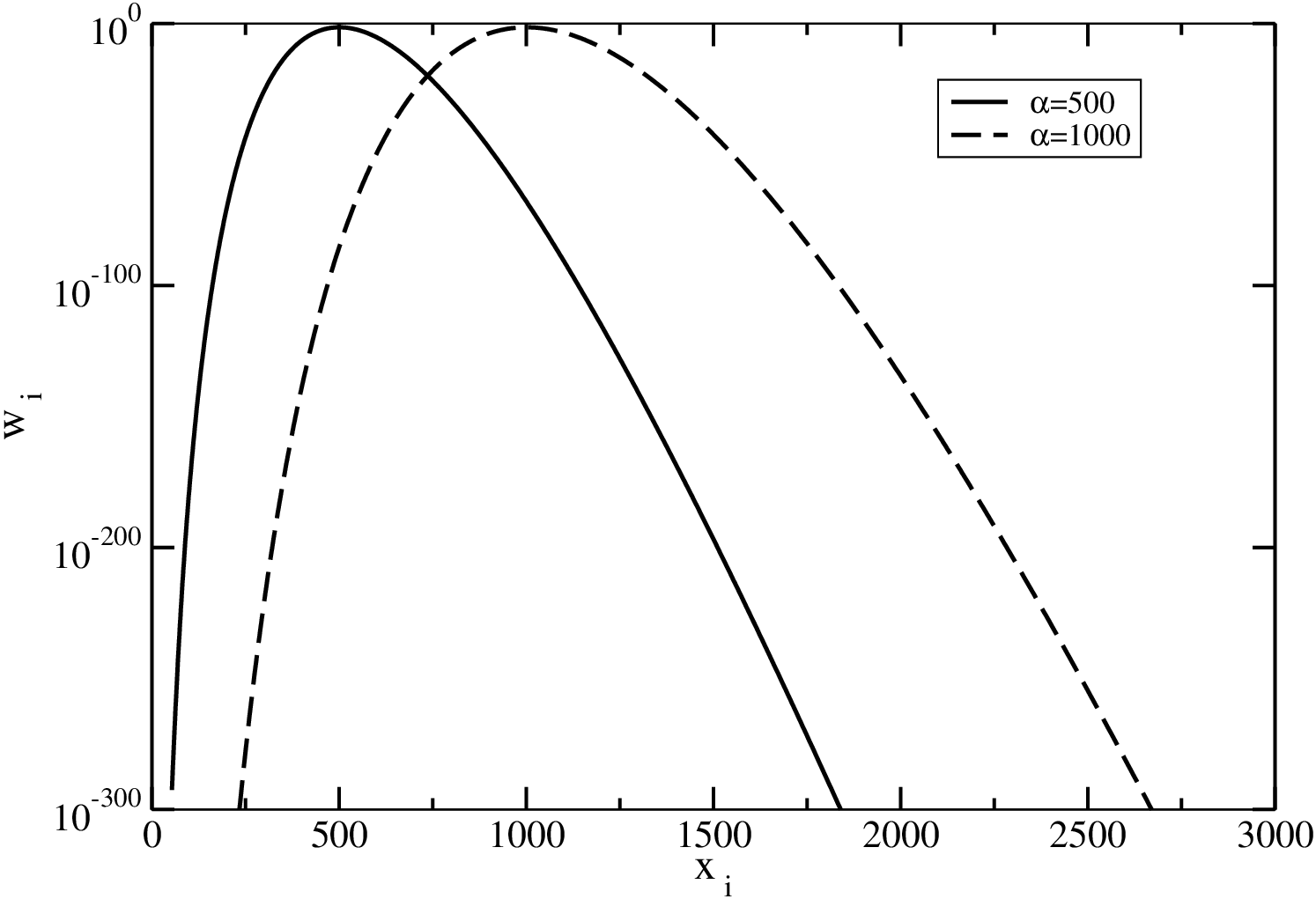}}
\caption{Gauss--Laguerre unscaled weights $\hat{w}_i$ as a function of the nodes $x_i$, 
where the degree is $n=1000$. Two values of $\alpha$ are 
considered: $\alpha=500$ and $\alpha=1000$}
\label{we1}
\end{figure}

However, as we explained before, the algorithm computes the weights in decreasing
order of significance, and the error is smaller for the first zeros and nodes
computed, larger or smaller than $x_e=\sqrt{\alpha^2 -1/4}$, when $|\alpha |>1/2$. 
The algorithm starts
at $x_e$, were the unscaled weights $\hat{w}_i$ are larger (and for $|\alpha|\le 1/2$, 
as explained earlier, the situation
is similar in that the weights are computed in decreasing order of magnitude).
This is shown in Figure~\ref{we1}.

\begin{figure}
\vspace*{0.8cm}
\centerline{\includegraphics[height=6cm,width=12cm]{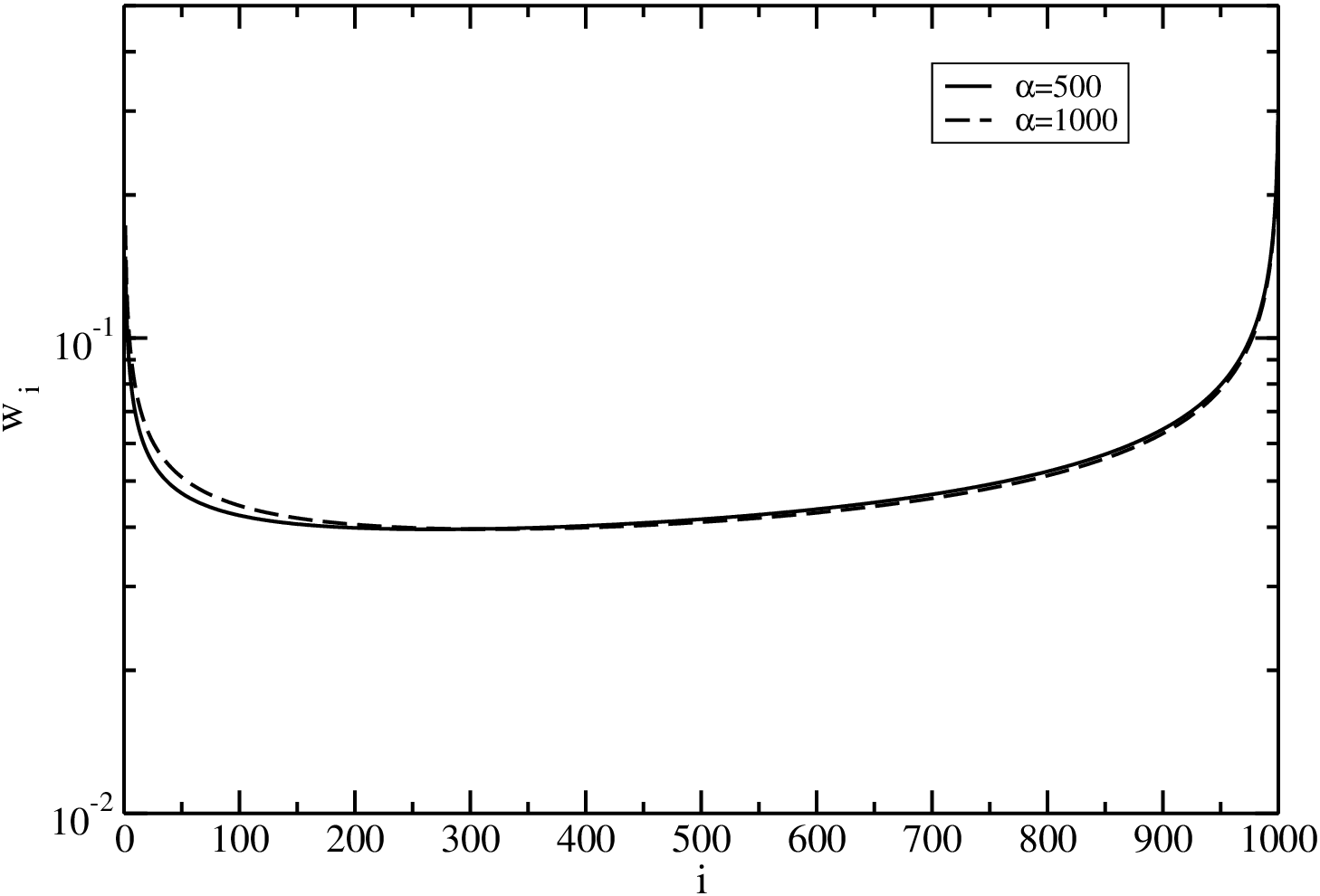}}
\caption{Gauss--Laguerre scaled weights $\hat{\omega}_i$ as a function of the nodes $i$ for
$n=1000$ and two values of $\alpha$: $\alpha=500$ and $\alpha=1000$}
\label{we2}
\end{figure}

The scaled weights, contrarily, and as expected, have a much smoother variation,
as they are approximately constant for large $n$. This is shown in Figure~\ref{we2}, 
were the scaled weights $\hat{w}_i$ are represented as a function
of $i$. We observe that the dependence on $\alpha$ is also very smooth and that
both curves are close to be indistinguishable.

\begin{figure}
\vspace*{0.8cm}
\centerline{\includegraphics[height=6cm,width=12cm]{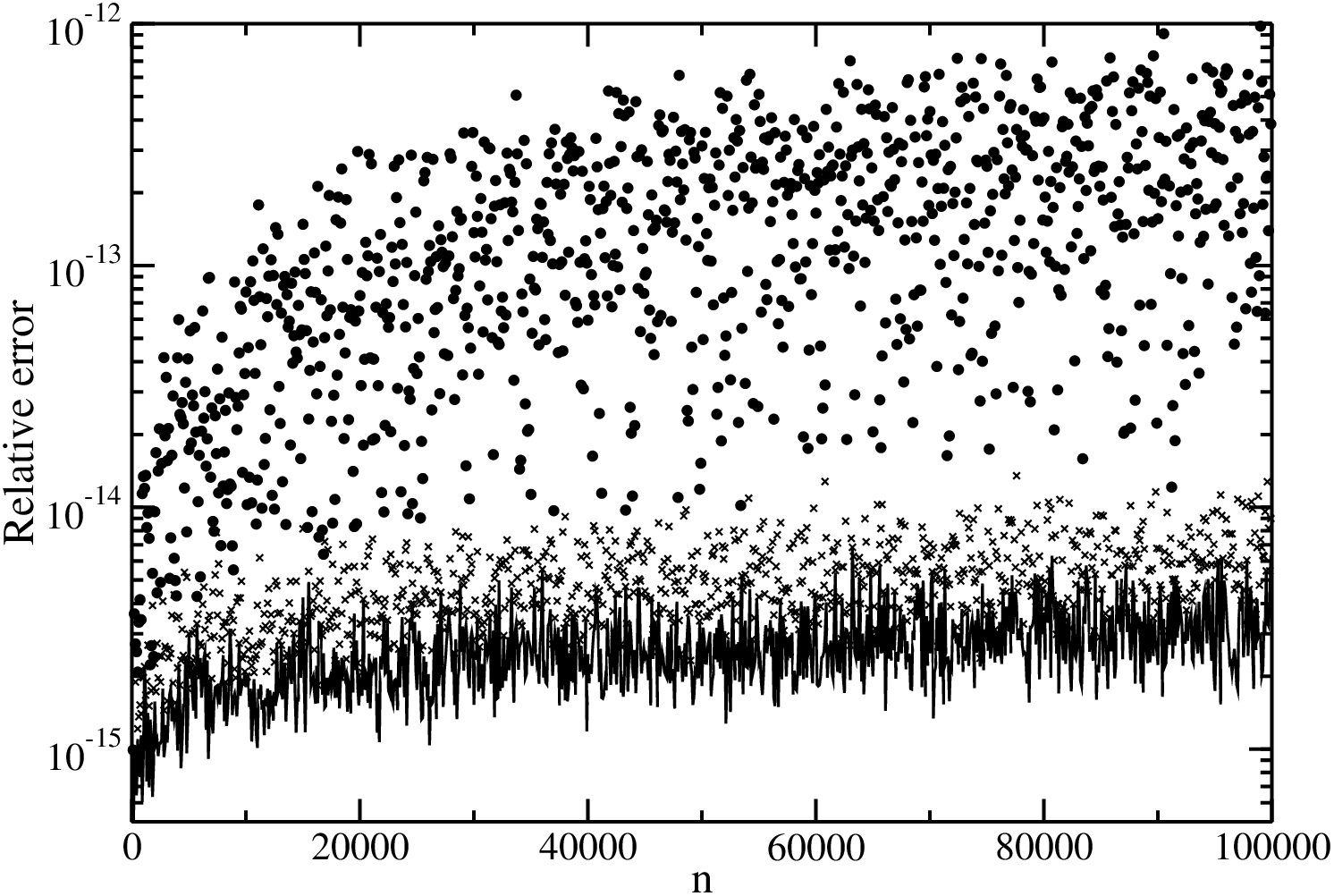}}
\caption{Relative errors in the computation of the nodes for $n$-point Gauss--Laguerre
quadrature with $\alpha=0$. The dots represent the values 
$\max|1-x_{i}^{(d)}/x_i^{(q)}|$, where $x_i^{(d)}$ are the nodes
 computed in double precision and 
$x_i^{(q)}$ are the same nodes in quadruple precision. 
We also show 
the maximum error when it is evaluated
 only for the nodes for which the (unscaled) weights
$w_i$ are larger than $10^{-300}$ (crosses) and $10^{-30}$ (solid line); the errors in these cases are smaller.}
\label{xnol}
\end{figure}

\begin{figure}
\vspace*{0.8cm}
\centerline{\includegraphics[height=6cm,width=12cm]{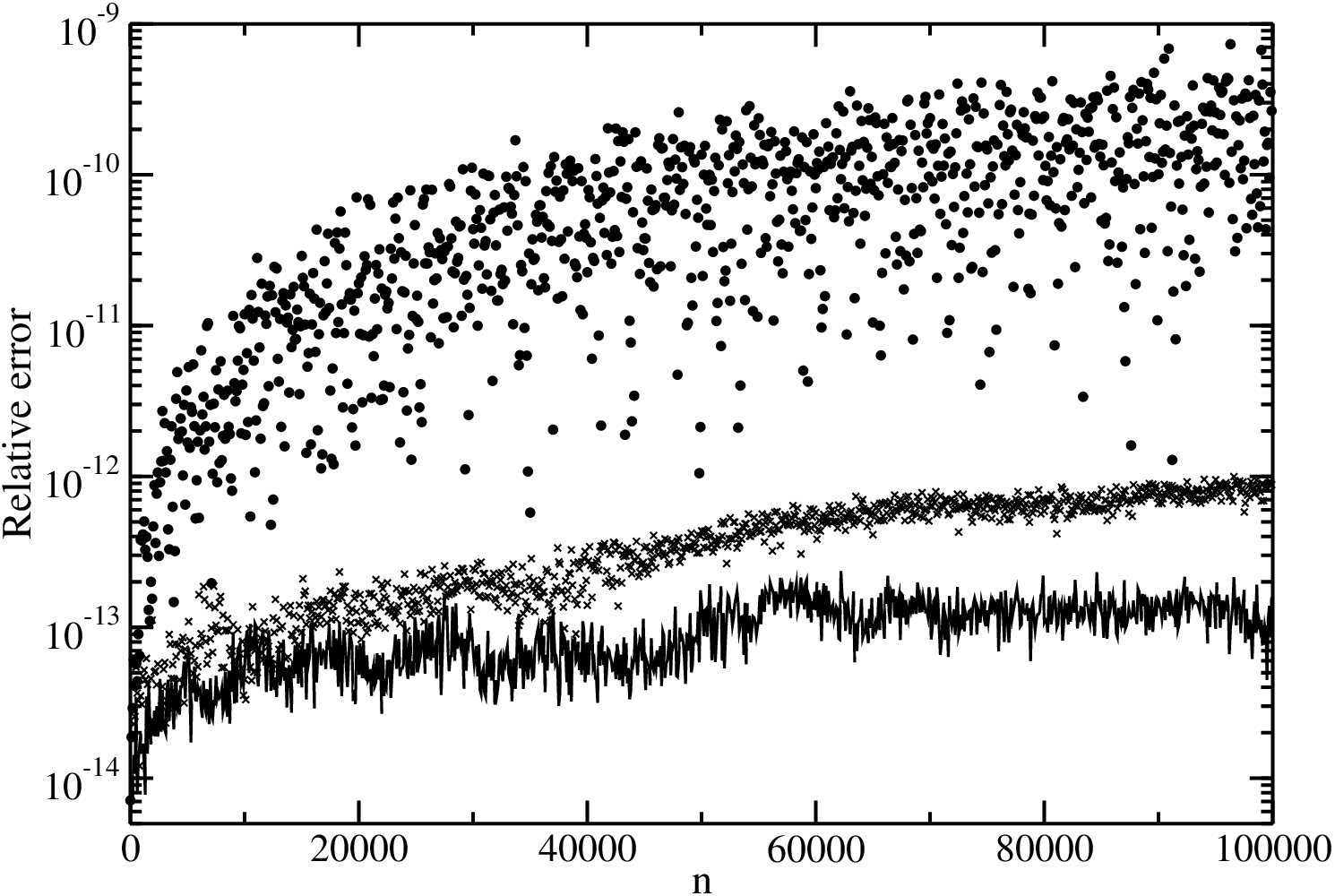}}
\caption{Relative errors in the computation of the scaled weights
 for $n$-point Gauss--Laguerre
quadrature with $\alpha=0$. The dots represent the values 
$\max|1-\hat{\omega}_{i}^{(d)}/\hat{\omega}_i^{(q)}|$, where 
$\hat{\omega}_i^{(d)}$ are the weights
 computed in double precision and 
$\hat{\omega}_i^{(q)}$ are the same weights in quadruple precision. 
We also show 
the maximum error when it is evaluated
 only when the (unscaled) weights
$w_i$ are larger than $10^{-300}$ (crosses) and $10^{-30}$ (solid line); in these cases the errors are smaller.}
\label{xwe}
\end{figure}

The error, for both the nodes and the weights is, as commented before, larger
as we move away from $x_e$ when $|\alpha|>1/2$ (or from $x=0$ for $|\alpha|\le 1/2$). 
Then, the error is larger as the weights become less
significant. This is shown in Figure~\ref{xnol}, where we plot the maximum relative
errors for the nodes as a function of $n$. We show three curves; one of them is the
maximum error considering all the weights, and the other two only considering
those nodes for which the unscaled weight $\tilde{w}_i$ are larger than $10^{-300}$ or $10^{-30}$.
We show the results for $\alpha=0$, but for other values of $\alpha$ the 
situation is similar (for instance, for $\alpha=100$ the results are almost
indistinguishable from those for $\alpha=0$). Figure~\ref{xwe} shows analogous
results, but for the unscaled weights $\hat{\omega}_i$.

\begin{figure}
\vspace*{0.8cm}
\centerline{\includegraphics[height=6cm,width=12cm]{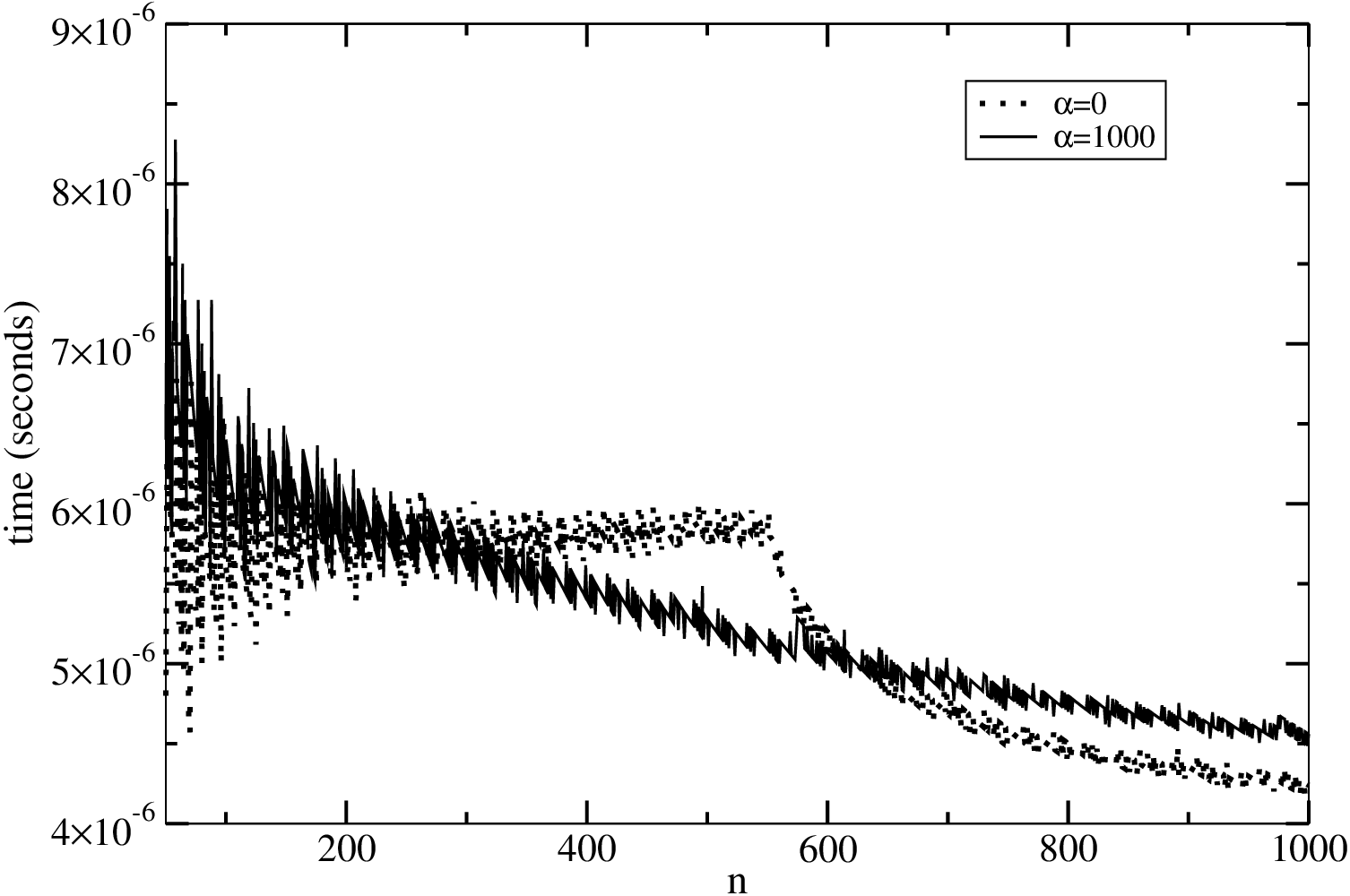}}
\caption{Unitary CPU-time spent as a function of the degree $n$ for Gauss--Laguerre with $\alpha=0$
and $\alpha=1000$}
\label{figCL}
\end{figure}

The algorithm for Laguerre quadrature is very efficient, but not so much as 
the one for the Hermite case (unsurprisingly). Figure~\ref{figCL}
shows the unitary time as a function of $n$ for two selections of the parameter 
$\alpha$.

Comparing with the CPU times spent by the asymptotic methods of \cite{Gil:2018:GHL}, we 
conclude that the asymptotic methods are faster by a factor smaller than $10$, and that
they are also more accurate. However, the present iterative method has several advantages 
with respect to the asymptotic methods. Firstly, the method is valid for any degree $n$, not 
necessarily large degree. Also, it is not limited to small $\alpha$, as are the methods considered
in \cite{Gil:2018:GHL}; in fact, this method is practically unrestricted with respect to $\alpha$,
which is a unique feature of the method. Finally, given that the
method is based on convergent approximations, it can be used for arbitrary accuracy (and we show some
results for very high accuracy in the next section); as an example of this, we point out 
that we have used a quadruple version of our algorithm to test our double precision implementation and
and that also the asymptotic methods in \cite{Gil:2018:GHL} have been tested against our iterative methods.
The same could be said with respect to other types of asymptotic methods, like for instance those based
in the Riemann--Hilbert approach of \cite{Huy:2018:CAI} (however those types of techniques can also be
considered for non-classical weights). With respect to the fully iterative method of 
\cite{Glaser:2007:AFA}, to the advantages
 already discussed also for the Hermite case (faster and certain higher order convergence and higher 
accuracy), we must add that our algorithm is not restricted to $\alpha=0$, and that in fact it works for
practically unrestricted $\alpha$.

\section{Perspectives and further applications}

In a next publication, we will consider the iterative computation of Gauss--Jacobi quadratures as well
as Gauss--Radau and Gauss--Lobatto quadratures.

Gauss--Jacobi quadrature can be treated in a similar manner. However, there exists a
number of characteristics that are different and which require further analysis. To start with,
the interval of integration is finite and the clustering of nodes for high degrees poses an 
additional stability problem. In addition, differently from the Hermite and Laguerre cases,
the canonical variable for which the method becomes asymptotically exact, is not suitable 
for computing Taylor series as we did before, because the relation with the original variable $x$ 
is $x=\cos\theta$, and therefore the derivatives with respect to $\theta$ do not satisfy a recurrence relation
with a fixed number of terms. In this case, it is likely that the computations will combine both the use
of the original and the canonical variable, and other possible changes of variables (particularly those described
in \cite{Deano:2004:NIF}). The initial values for starting the computation will also be necessarily more involved
than in the Laguerre case, because we are dealing with an additional parameter. We postpone the analysis to a future
publication. Asymptotic approximations for the nodes and weights are discussed in a recent paper 
\cite{Gil:2018:GHL} and, as for the
Hermite and Laguerre cases, these estimations can be considered as a standalone 
alternative to the iterative method for high enough orders (provided the zeros of Bessel
functions, which are used in the asymptotic expansions, are available).

Gauss--Radau and Gauss--Lobatto quadratures can also be also computed by following similar schemes for the internal
nodes and computing the boundary nodes with the particular formulas for these cases (Gauss--Lobatto quadratures
do not make sense in the present case, because we are not dealing with finite intervals, but they will be
considered for the Gauss--Jacobi case). The generalized
Gauss--Radau--Laguerre quadrature formula is the approximation
$$
\displaystyle\int_{0}^{+\infty} f(x) x^{\alpha} e^{-x} dx \approx \displaystyle\sum_{j=0}^{r-1} 
w^{(j)}_0 f^{(j)}(0) +
\displaystyle\sum_{i=1}^n 
w^{(R)}_i f(x_i) ,
$$
with the highest possible degree of exactness, which is $2n-1+r$. As is well known, 
the internal nodes $x_i$ are the zeros of $L_n^{(\alpha+r)}(x)$ and the weights $w^{(R)}_i$ can be written 
in terms of the Gauss--Laguerre weights $w_i$ of degree $n$ and parameter $\alpha+r$ as 
$w^{(R)}_i =w_i /(x_i)^r$ \cite{Gautschi:2004:GGR}, while the boundary nodes $w^{(j)}_0$ can be 
computed using the methods in \cite{Gautschi:2004:GGR}. In particular, when $r=1$, we are dealing with 
the Gauss--Radau--Laguerre formula and the boundary weight is explicitly given 
(see \cite{Gautschi:2000:GRF}) by
$$
w_0^{(0)}=\Gamma (\alpha +1)\left/
\left(
\begin{array}{c}
n+\alpha+1\\
n
\end{array}
\right).
\right.
$$
Therefore, the algorithms we have constructed in this paper are also of application to (generalized) 
Gauss--Radau--Laguerre quadrature because the internal nodes and weights can be computed using the same scheme.

In addition, as pointed out in \cite{Hale:2013:FAA}, the computation of Gauss quadrature rules is related to
the problem of interpolation at the orthogonal polynomial nodes with the barycentric formula. Indeed, the
Lagrange interpolation polynomial at the simple zeros $x_i$ of a polynomial $q(x)$ of degree $n$ 
for a function $f(x)$ can be written
\begin{equation}
\label{bar}
P_{n-1}(x)=\displaystyle\sum_{i=1}^{n}\Frac{v_i f_i}{x-x_i}\left/
\displaystyle\sum_{i=1}^{n}\Frac{v_i}{x-x_i}\right.
\end{equation}
and this is the lowest degree polynomial (of degree $n-1$) satisfying the interpolating conditions 
$P_{n-1}(x_i)=f(x_i)=f_i$ when the weights $v_i$ are computed by $v_i=1/q'(x_i)$ (an additional constant 
factor for all weights can be also considered). We observe that our algorithms allow us
to interpolate functions at the Hermite and Laguerre nodes even for very high degrees and for practically 
unrestricted values of $\alpha$ for the Laguerre case. Indeed, we compute both the nodes $x_i$ and 
the derivative of the polynomial $q(x)$ ($H_n(x)$ and $L_n^{(\alpha)}(x)$ in our case), up to an elementary 
scale factor, say $s(x)$: we therefore can use the derivative of $y(x)=s(x)q(x)$ (solution of the second
order ODE in normal form) to compute $q'(x_i)$. Notice that, as commented, the algorithm computes Gaussian weights in the direction of decreasing weights, which is also the direction of decreasing values of $v_i$. This,
 together with the possible use of scaling factors for the function to be interpolated, 
is an interesting property in order to avoid underflows
in the evaluation of (\ref{bar}).

\section{Conclusions}

We have described fast and reliable iterative methods for the 
computation of Gauss--Hermite and Gauss--Laguerre quadratures. These methods have a number of 
interesting and distinctive features, among them:

\begin{enumerate}

\item The computation of the nodes is based on a globally convergent fourth-order method. The
convergence is certain and fast. No initial estimations for the nodes are needed. 

\item The methods are valid for small and large degrees.

\item Choosing what we called the canonical variable, the iterative method is asymptotically exact
as the degree goes to infinity, and the computational time per node decreases as the degree increases.

\item In the canonical variable, we have defined well-conditioned scaled weights.

\item The methods are essentially unrestricted with respect to the range
of the parameters, thanks to weight scaling and normalization of solutions of the ODE.

\item The weights are computed in decreasing order of magnitude and the most significant weights are the
most accurate ones. This is useful for subsampling (computing only nodes and weights
corresponding to weights greater than a given threshold).

\item Because the methods only use convergent procedures, they can be used for arbitrary accuracy. The
fast fourth-order convergence makes this an interesting method for high accuracy computations.

\end{enumerate}

\begin{acknowledgements}
The authors thank the anonymous referees for their constructive comments and suggestions. NMT thanks CWI  for scientific support.
\end{acknowledgements}

\bibliographystyle{spmpsci}      
\bibliography{gaussr}   

\end{document}